\documentclass[dvips,11pt,letter]{article}

\usepackage{amssymb,epsfig,latexsym}
\usepackage{graphicx,amsmath,amssymb,latexsym}
\usepackage{setspace}
\def\cV{{\cal V}}
\def\cA{{\cal A}}
\def\cX{{\cal X}}

\def\cK{{\cal K}}
\def\cL{{\cal L}}

\def\cY{{\cal Y}}
\def\cZ{{\cal Z}}
\def\cC{{\cal C}}
\def\cR{{\cal R}}

\def\st{{\rm st}}
\def\cF{{\cal F}}

\def\N{\mathbb{N}}
\def\R{\mathbb{R}}

\def\cB{{\cal B}}
\def\cP{{\cal P}}
\def\cS{{\cal S}}
\def\and{\quad\mbox{and}\quad}
\def\ind{{\bf 1}}

\def\PP{{\mathrm P}}
\def\EE{{\mathrm E}}

\newtheorem{theorem}{Theorem}[section]

\newtheorem{coro}[theorem]{Corollary}
\newtheorem{lemma}[theorem]{Lemma}

\def\bp{\noindent{\it Proof. }}
\def\ep{\hfill $\Box$}

\begin{document}

\title{A particle system in interaction with a rapidly varying environment: Mean field limits and applications}

\author{Charles Bordenave\footnote{CNRS, Universit\'e de Toulouse, France.} ,
David McDonald\footnote{University of Ottawa, Canada.} \,  and Alexandre Proutiere\footnote{Microsoft Research, Cambridge, UK.}}

\maketitle

\begin{abstract}
We study an interacting particle system whose dynamics depends on an
interacting random environment. As the number of particles grows
large, the  transition rate of the particles slows down (perhaps
because they share a common resource of fixed capacity). The
transition rate of a particle is determined by its state, by the
empirical distribution of all the particles and by a rapidly varying
environment. The transitions of the environment are determined by
the empirical distribution of the particles. We prove the
propagation of chaos on the path space of the particles and
establish that the limiting trajectory of the empirical measure of
the states of the particles satisfies a deterministic differential
equation. This deterministic differential equation involves the time
averages of the environment process.

We apply the results on particle systems to understand the behavior of computer networks where users access a shared resource using some distributed random Medium Access Control (MAC) algorithms. These algorithms are used in all Local Area Network (LAN), and have been notoriously difficult to analyze. Our analysis allows us to provide, for the first time, simple and explicit expressions of the network performance under such algorithms.



\end{abstract}

\section{Introduction and motivation}

The paper comprises two separate parts: a first part is devoted to the analysis of the mean field limits of a general system of interacting particles, also interacting with a random environment. In the second part of the paper, we demonstrate how the results on particle systems derived in the first part can be applied to understand the behavior of computer networks where users access a shared resource using some distributed random Medium Access Control (MAC) algorithms. These algorithms are implemented in the network access card of all computers connected to a Local Area Network (LAN). LANs are networks covering a small geographic area, like a home, an office, a building, and constitute a first and crucial component of the Internet. Analyzing random MAC algorithms is notoriously difficult; most of the related issues have actually been open since the introduction of the first of these algorithms in the early 70's. In \cite{bianchi}, the author has used heuristic formulas to approximate their performance in specific networks. These formulas are based on the assumption that  the particles (or computers) evolve independently. Our  mean field analysis rigorously proves this propagation of chaos, but also allows for the first time to derive explicit analytical expressions of the performance of these algorithms in general networks.

\paragraph{A particle system interacting with a random environment}
In the first part of the paper, we are  interested in the mean field limit of a system of $N$ interacting particles whose dynamics also depends on an environment process. More specifically,  the evolution of each particle depends on the state of the particle, on the empirical distribution of all the particles and also on environment variables. The environment process is a finite state space Markov chain which interacts with
the particle system because its transition kernel  depends on the empirical distribution of the states of the particles.  A key feature of the systems considered here is that the  environment is rapidly varying: it evolves at rate $1$, whereas the particles evolve at rate $1/N$.

We prove a mean field limit for this particle system when $N$
goes to infinity.  In order to capture the evolution of the particles we must speed up time by a factor of $N$.
In so doing the particles see a time average of the rapidly changing environment. In the mean field limit, particles evolve independently, and see the environment process in its steady state, which in turn evolves as the particles evolve. 

 Our  results on a particle system evolving in a rapidly changing environment
is a  generalization of results obtained by Kurtz in \cite{kurtzaverage}. We extend these results in a couple of ways: first, the particles may evolve according to their current states, to the empirical distribution of the states of the particles, and to an environment process; then we show the path-space convergence of the trajectory of the empirical distribution of the states of the particles. To prove this convergence, we extend and adapt the method developed by Sznitman and Graham in \cite{sznitman,graham}. 
%In general deriving conditions for global stability of the mean field limit is a hard problem.  Here, the best we could do (still, no small achievement) was to prove global stability of the mean field limit
%of a system of Wi-Fi users in total interaction (see Section \ref{globalstab}).

The initial motivation for the use of mean field asymptotics was to analyze the behavior of computer networks. Of course, mean field models have been used in many contexts and the theory is well developed. For example,  Dawson \cite{dawson} studies a model in statistical physics where  $N$ particles diffusing in a potential well have the additional property that they are all attracted to the center of mass of all the particules. The Fleming-Viot model \cite{donnelly} is an example from genetics where a particle represents an individual and its state represents the genetic type and its location. Our results (and those in \cite{kurtzaverage}) on  interacting particle systems with a rapidly varying environment could find other applications.   For instance, they could be used to capture the dynamics of a population whose genetic makeup evolves slowly in time in the presence of a rapidly varying environment whose evolution may partly depend on the empirical distribution of the individuals. Another potential field of application is microscopic models in
economic theory and stochastic market evolution, also known as
"econophysics", see for example the work by Karatzas \cite{karatzas}
or Cordier \cite{cordier}. In a simple market economy or in a
financial market, a particle is an economic agent and its states
represents its goods and its savings. The environment is the prices
of the various available goods. Agents may exchange, borrow or lend
money. Both prices and the purchase decisions of agents are
interacting. In some markets, like financial markets, the prices are
fluctuating roughly $N$ times faster than the decisions of each
individual agent.

\paragraph{Analyzing Medium Access Control algorithms in computer networks}

Consider  $N$ users (or computers) communicating in a wired or wireless Local Area Network (LAN). To transmit data packets, users have to share a single resource (a cable in wired LANs or a radio channel in wireless LANs) using some Medium Access Control (MAC) protocols. These protocols are distributed, meaning that each user runs its protocol independently of the other users sharing the same resource. This architecture has ensured the {\it scalability} of LANs (in the sense that new users can join and leave the network without the need of explicitly advertising it); it has played a crucial role in their development and hence contributed to the rapid growth of the Internet.

When two users cannot simultaneously successfully transmit data packets (because they share the same resource), we say that these users {\it interfere}. Two interfering users who simultaneously transmit experience a {\it collision}, and the packets have to be retransmitted. Most current MAC protocols  limit collisions using the following two main principles: first, before transmitting, users sense the resource and should it be busy they abstain from transmitting. This technique is referred to as CSMA (Carrier Sense Multiple Access) and ensures that packet transmissions cannot be interrupted. Even if the sensing mechanism is perfect, a collision may still occur if two interfering users start transmitting at the same time (or rather so close together in time that CSMA can't prevent the collision). The second main principle, termed random back-off, aims at reducing the possibilities that several users start transmitting simultaneously. To do so, a user only starts  transmitting with a certain probability less than one. This probability is adapted to the number of successive collisions experienced by users, which allows users to infer the level of congestion of the resource. Typically, in LANs today, users implement the exponential back-off algorithm (also referred to as the Decentralized Coordination Function (DCF) in the standards, see \cite{bianchi} and references therein for a detailed description of these standards): the transmission probability is divided by a factor two after each collision, and it is reinitialized after the successful transmission of a packet.

The performance of MAC protocols is measured in terms of the throughput realized by the various users, i.e., of the number of packets successfully transmitted by users per second. The performance analysis requires that we can characterize the joint evolution of the transmission probabilities of the $N$ users (see Section \ref{sec:example} for the state of the art). These probabilities evolve according to a $N$-dimensional Markov chain that is usually intractable because of the correlations introduced by collisions. Mean field asymptotics are useful to approximate this evolution.  

\begin{figure}[h]\label{firstpic}
\centering
\includegraphics[width=8cm]{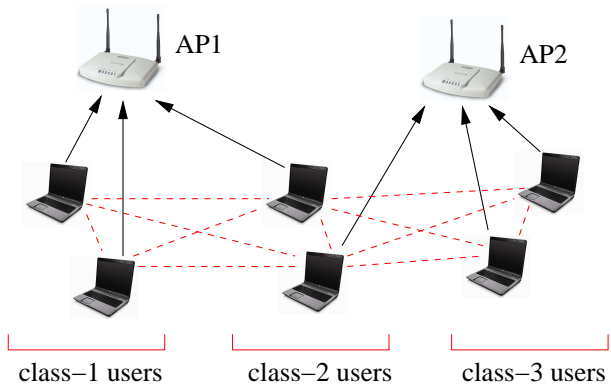}
\caption{A network with partial interference - A dashed line between two users means that they sense each other.}
\end{figure}

In this paper, we consider two relevant scenarios for interference. We consider networks with {\it full interference} where all pairs of users interfere, and networks with {\it partial interference} where users do not interfere with all other users. In the latter scenario, users are classified according to the set of users they interfere with. Partial interference typically arises in wireless networks as illustrated in Figure \ref{firstpic}: all 6 users are willing to transmit data packets to the access points 1 or 2; class-1 (resp. class-3) users interfere with users of classes 1 and 2 (resp. 2 and 3), whereas class-2 users interfere with all users. Two users of class 1 and 3 respectively can not sense each other and this can lead to fairness issues: users of class 2 find themselves in a  predicament like that of a polite nephew sitting on a sofa between two garrulous aunts who are hard of hearing and therefore hear the nephew but not each other. Each aunt will launch into a new dialogue before the other aunt has finished.  The poor nephew will hardly ever get a word in!

In Section \ref{sec:example}, we apply the results derived for the particle systems to provide accurate approximations of the performance of a general class of MAC protocols in networks with full or partial interference. A user is modeled as a particle whose state includes its transmission probability and its class in the case of partial interference. The particles interact with each other because of collisions. If the access protocol is fair each user will necessarily share  around $1/N$ of the resource; i.e. the transmission probability of users slows down as $N$ increases. Our mean field limit  will therefore depend on rescaling and speeding up time by a factor of $N$. The environment process captures the fact that for a given user, the resource is sensed busy or idle. For example, the environment of the network in Figure \ref{firstpic} is represented by a vector $z=(z_1,z_2,z_3)$ of zeros and ones.  The environment $(1,0,1)$ would represent ongoing transmissions from a user in class 1 and a user in class 3. When a user transmits the resource is blocked; i.e. the environment changes.  These environmental changes occur at rate determined by $N$ users; i.e. at rate $1$. Consequently the conditions for our theory are met.

\paragraph{Notations}  Let $S$ be a separable, complete metric space, ${\cal
P}(S)$ denotes the space of probability measures on ${\cal
Y}$. ${\cal L}(X)$ is the law of the $S$-valued random
variable $X$. $D(\R^+,S)$ the space of right-continuous
functions with left-handed limits, with the Skorohod topology
associated with its usual metric, see \cite{ethier} p
117. With this metric, $D(\R^+,S)$ is complete and separable.
We extend a discrete time trajectory $(X(k)), k \in \N,$ in
$D(\N,S)$ in a continuous time trajectory in $D(\R^+,S)$ by setting for $t \in \R^+$, $X(t) = X ([t])$, where $[\cdot]$
denotes the integer part. $(\cF_t), t \in \R^+ \hbox{ or } \N,$ will
denote the natural filtration with respect to the processes
considered. $\| \cdot\|$ denotes the norm in total variation of
measures. Finally, for any measure $Q\in {\cP}(S)$ and any
measurable function $f$ on $\cY$, $\langle f, Q \rangle = Q(f) = \int f dQ$ denotes the usual duality brackets.

We recall that a sequence of random variables
$(X_i^N)_{i\in\{1,\ldots,N\}}\in S^N$ is exchangeable if
$\cL((X_i^N)_{i\in
\{1,\ldots,N\}})=\cL((X_{\sigma(i)}^N)_{i\in\{1,\ldots,N\}})$ where
$\sigma$ is any permutation of $\{1,\ldots,N\}$. Moreover the
sequence is $Q$-chaotic if for all subsets $I\subset \N$ of finite
cardinal $\vert I\vert$,
\begin{equation}\label{q-chaotic}
\lim_{N\to \infty} {\cal L}\left((X_i^N)_{i\in I}\right)= Q^{\otimes
\vert I\vert} \quad \hbox{weakly in } {\cal P}(S^{\vert I\vert}).
\end{equation}

%$C$ will denote a constant not depending on the parameters of
%interest of the system.

\section{An interacting particle system in a varying
environment}\label{markov}

In this section, we first provide a precise description of the
interacting particle system under consideration. We then state the
main results, giving the system behavior in the mean field limit
when the number of particles grows to infinity. The proofs of these
results are postponed to subsequent sections.

\subsection{Model description}
\label{subsec:model}

\paragraph{The particles} We consider $N$ particles evolving in a countable state space
$\cX$ at discrete time slots $k\in\mathbb N$. For simplicity we
assume the particles are exchangeable. At time $k$, the state of the
$i$-th particle is $X_i ^N (k) \in \cX$. The state of the system at
time $k$ is described by the empirical measure $\nu^N(k)\in
\cP(\cX)$ while the entire history of the process is described by
the empirical measure $\nu^N$ on path space $\cP(D(\mathbb N,\cX))$:
$$\nu^N(k) = \frac 1 N \sum_{i=1} ^ N \delta_{X^N_i(k)}\quad\mbox{ and
}\quad \nu^N = \frac 1 N \sum_{i=1} ^ N \delta_{X^N_i}.$$

\paragraph{The interacting environment} In the system considered, the evolution of the
particles depends not only on the state of the particle system but
also on a background Markovian process $Z^N = (Z^N_1,\cdots,Z^N_N)
\in D(\mathbb N,\cZ^N)$, where $\cZ$ is a countable state
space. Specifically, $Z^N$ is a Markov chain whose transition kernel satisfies the following:
\begin{equation}\label{eq:ker}
\PP ( Z_i^N(k+1)=z | \cF_k ) = K^N_{\nu^N(k),X^N_i(k)}
(Z^N_i(k),z),
\end{equation}
where $K^N_{\mu,x}$ is a transition kernel on $\cZ$ depending on a
probability measure $\mu$ on $\cP(\cX)$ and on $x\in \cX$, and where $\cF_k = \sigma
\left((\nu^N(0),Z^N (0)),\cdots,(\nu^N(k),Z^N (k)) \right)$. The
latter filtration depends on $N$, but as pointed out above,
without possible confusion,  $\cF_k$ will always denote the
underlying natural filtration of the processes. Note that (\ref{eq:ker}) does not completly defined the transition kernel of $Z^N$, and actually the joint
evolution of the vector $(Z^N_1 (k),\cdots,Z^N_N(k))$ is arbitrary.

\paragraph{Evolution of the particles} We represent the possible
transitions for a particle by a countable set $\cS$ of mappings
from $\cX$ to $\cX$. A $s$-transition for a particle in state $x$
leads this particle to the state $s(x)$. We assume that the
conditional probability given $\cF_k$ that a $s$-transition occurs
for the particle $i$ between times $k$ and $k+1$ is equal to
\begin{eqnarray}\label{jumps}
\frac 1 N F^N _s (X^N_i(k),\nu^N(k),Z^N_i (k)).
\end{eqnarray}
with $\sum_{s\in\cS} F^N_s(x,\alpha,z) = 1$ for all $(x,\alpha,z)
\in \cX \times \cP(\cX) \times \cZ$ (the assumption is for
simplicity, the content of the paper is unchanged if
$\sum_{s\in\cS} F^N_s(x,\alpha,z) \leq C$ for some constant $C$
independent of $(x,\alpha,z)$).

We define the events $$A^N_i (k)= \{ \hbox{a transition occurs for
particle $i$ between times $k$ and $k+1$}\}.$$ We assume that the
joint distribution of the transitions is weakly correlated. More
precisely, 
\begin{itemize}
\item[A0.]
There exists a positive sequence $(\rho_N)_{N \in \N}$
such that $\lim_N \rho_N = 0$ and
\begin{equation}
\label{cancellation} \PP(A^N_1 (k) A^N_2(k) | \cF_k )   \leq \frac
{\rho_N}{N}.
\end{equation}
\end{itemize}

Note that, due to (\ref{jumps}), the process $Z^N$ evolves quickly
while the empirical measure $\nu^N(k)$ evolves slowly. Also note
that the $s$-transitions of the various particles may be
correlated. The process $Z^N$ may depend on the transitions of the
particles. The particle system is thus in interaction with its
environment. Note finally that if the particle transitions are
independent then (\ref{cancellation}) holds with $\rho_N =1/N$.

We make the following additional assumptions on the system
evolution.

\paragraph{Assumptions}
\begin{itemize}
\item[A1.] Uniform convergence of $F_s^N$ to $F_s$:\\
$\lim_{N \to \infty} \sup_{(x,\alpha,z) \in \cX\times\cP(\cX)\times
\cZ } \sum_{s \in \cS} |F^N_s (x,\alpha,z) - F_s (x,\alpha,z)| = 0.$
\item[A2.] The functions $F_s$ is uniformly Lipschitz:\\
$\sup_{(x,z) \in \cX\times \cZ } \sum_{s \in \cS} |F_s ( x, \alpha,
z) - F_s ( x, \beta, z) | \leq C \| \alpha - \beta\|$.
\item[A3.] Uniform convergence in total variation of $K_\alpha^N$ to
$K_\alpha$:\\
$\lim_{N \to \infty} \sup_{(x,\alpha,z) \in \cX\times
\cP(\cX)\times \cZ } \|K^N_{\alpha,x} (z, \cdot) - K_{\alpha,x}
(z, \cdot)\| = 0$.
\item[A4.] The mapping $\alpha \mapsto K_{\alpha}$ is uniformly Lipschitz: \\
$\sup_{(x,z) \in  \cX \times \cZ }  \|K_{\alpha,x} ( z, \cdot ) -
K_{\beta,x} ( z, \cdot) \| \leq C \| \alpha - \beta\|$. 
\item[A5.]
The Markov chains with kernels $K_{\alpha,x}$ have a unique
stationary probability measure $\pi_{\alpha,x}$. 

\item[A6.] For all $x$ in $\cX$, $\alpha$, $\beta$ in $\cP(\cX)$:
$\| \pi_{\alpha,x} - \pi_{\beta,x} \| \leq C \sup_{z \in  \cZ }
\|K_{\alpha,x} ( z, \cdot ) - K_{\beta,x} ( z, \cdot) \|.$

\end{itemize}

We discuss in Section \ref{sec:HYP} how the above assumptions may be
checked.

\subsection{Main Results}\label{subsec:results}

The main result of this paper is to provide a mean field analysis of
the system described above, i.e, to characterize the evolution of
the system when the number of particles grows. According to
(\ref{jumps}), as $N\to\infty$, the chains $X^N_i(t)$ slow down
hence to derive a limiting behavior we define:
$$
q_i ^ N (t) = X^N_i([Nt])\quad \hbox{ and }\quad \mu^N = \frac 1 N
\sum_{i=1} ^ N \delta_{q^N_i} \in \cP(D(\mathbb R_+,\cX)).
$$

We wish to apply the ideas in Theorem 2.1 in \cite{kurtzaverage}. In that context
we define the joint measure 
$$\zeta^N(k)(A\times B)=\frac{1}{N}\sum_{i=1}^{N}\chi\{X^N_i(k)\in A,Z^N_i(k)\in B\}$$
for $A \subseteq \mathcal{X}$ and $B \subseteq \mathcal{Z}$.
Clearly the evolution of $\nu^N$ is determined by $\zeta^N$. Next we rescale time and define
$Y^N(t)(A\times B)=\zeta^N([Nt])$. In the context of \cite{kurtzaverage} our $\nu^N$ is Kurtz's $X_N$
and our $Y^N$ is Kurtz's $Y_N$. However we can't quite apply the theorems in \cite{kurtzaverage}
because the transition kernel of  $Z^N_i$ depends on both $\nu^N$ and $X^N_i$.

Following \cite{kurtzaverage} we define $\ell_m(\cZ\times\cX)$ to be the space of
measures on $[0,\infty)\times\cZ\times\cX$ such that for $\gamma\in \ell_m(\cZ\times\cX)$,
$\gamma([0,t]\times\cZ\times\cX)=t$. Define 
$$\Gamma^N([0,t]\times A\times B)=\int_0^tY^N(s)(A\times B)ds.$$
Note that $\Gamma^N([0,T],x,\cZ)=\int_0^T\mu^N(s)(x)ds.$
Since $Y^N$ doesn't slow down as $N\to\infty$ like $\mu^N$ we can't hope to prove the weak
convergence of $Y^N$ but
 the occupation measure $\Gamma^N$ does converge weakly by averaging.  To
obtain the relative compactness of $\Gamma^N$ and $\mu^N$ we require the following assumptions.

\begin{itemize}
\item[A7.] For each $\epsilon>0$ and each $t>0$ there exists a compact $\cK\subseteq \cX\times\cZ$
such that $\liminf_N E[\Gamma^N([0,t]\times \cK)]\geq (1-\epsilon)t$.
\item[A8.]  $\cL(q_1^ N (\cdot))$ is tight in $\cP(D(\R^+,\cX))$.
\end{itemize}

In most applications, the tightness of $\cL(q_1^N (\cdot))$ in $\cP
(D(\R^+,\cX))$ is not a major issue. Indeed, note that the
inter-arrival times between two transitions of $q_1^N(.)$ are
independent Binomial $(N,1/N)$ variables (which converges to
exponential $(1)$ variables).  Hence, if for example the state space
$\cX$ or the set of transitions $\cS$ is finite, we may apply the
tightness criterion Theorem 7.2 in Ethier-Kurtz \cite{ethier} p.128.

\subsubsection{Transient regimes}

The following theorem provides the limiting behavior of the system
in transient regimes.

\begin{theorem}\label{main}
Assume that the Assumptions A0-A8 hold and that the initial values
$q_i^N(0)$, $i=1,\ldots,N$, are exchangeable and such that their
empirical measure $\mu^N_0$ converges in distribution to a
deterministic limit $Q_0\in \cP(\cX)$ when $N\to\infty$. There exists a probability
measure $Q$ on $D(\R^+,\cX)$ such that the processes $(q_i^N(.),i\in \{1,\ldots,N\})$ are $Q$-chaotic.
\end{theorem}

In \cite{sznitman}, Sznitman proved that if $q_i^N(0)$, $i=1,\ldots,N$, are exchangeable, 
their empirical measure $\mu^N_0$ converges in distribution to a deterministic limit $Q_0\in \cP(\cX)$ 
if and only if $q_i^N(0)$, $i=1,\ldots,N$, are $Q_0$-chaotic. Then, the above
theorem states that if the particles are initially asymptotically
independent, then they remain asymptotically independent. This
phenomenon is also known as the propagation of chaos.

The independence allows us to derive an explicit expression for the
system state evolution. As explained earlier, intuitively, when $N$
is large, the evolution of the background process is very fast
compared to that of the particle system. The particles then see a
time average of the background process. The following theorem
formalizes this observation. For $\alpha \in \cP(\cX)$ and $x\in \cX$, let
$\pi_{\alpha,x}$ denote the stationary distribution of the Markov
chain with transition kernel $K_{\alpha,x}$. We define the average
transition rates for a particle in state $x$ by
\begin{equation}
\label{eq:Fbarre} \overline F_s (x,\alpha)= \sum_{z\in \cZ} F_s
(x,\alpha, z) \pi_{\alpha,x}(z).
\end{equation}

Define $Q^n(t) = Q(t)(\{x_n\})$ where $\cX = \{ x_n, n\in\N\}$.
$Q^n(t)$ is the limiting (when $N\to\infty$) proportion of particles
in state $x_n$ at time $t$.

\begin{theorem}\label{theoex}
Under the assumptions of Theorem \ref{main}, the limiting
proportions $Q^n(t)$ of the particles in the various states satisfy:
$Q^n(0)=Q_0(\{x_n\})$ and for all time $t>0$, for all $n\in \N$,
\begin{equation}\label{eqdiff1}
{d{Q^n}\over dt}= \sum_{s \in \cS} \sum_{m :
s(x_m) = x_n } Q^m (t) \overline F_s ( x_m, Q(t))  -\sum_{s \in \cS} Q^n (t)
\overline F_s ( x_n, Q(t)).
\end{equation}
\end{theorem}

The equations (\ref{eqdiff1}) have the following
interpretation: if $s(x_m) = x_n$ then $Q^m (t) \overline F_s ( x_m,
Q(t))$ is a mean flow of particles from state $x_m$ to $x_n$.
Hence, $\sum_{s \in \cS} \sum_{m : s(x_m) = x_n } Q^m (t) \overline
F_s ( x_m, Q(t)),$ is the total mean incoming flow of particle to
$x_n$ and $ \sum_{s \in \cS} Q^n (t) \overline F_s ( x_n, Q(t))$ is
the mean outgoing flow from $x_n$.

\subsection{Stationary regime}

We now characterize the stationary behavior of the system in the
mean field limit. To do so, we make two additional assumptions:
\begin{itemize}

\item[A9.] For all $N$, the Markov chain $((X_i^N(k))_{1 \leq i
\leq N},Z^N(k))_{k\in \N}$ is positive recurrent. The set of stationary distributions ${\cal L}_{st} (X^N_1)$ is tight.

\item[A10.] The dynamical system (\ref{eqdiff1}) is globally
stable: there exists a measure $Q_{\st} = (Q^n_{\st}) \in \cP(\cX)$
satisfying for all $n$:
\begin{equation}\label{eq:balance}
\sum_{s \in \cS}  \sum_{m : s(x_m) = x_n } Q^m _ {\st}\overline F_s (
x_m, Q_{\st}) =  Q^n _ {\st}\sum_{s \in \cS} \overline F_s ( x_n,
Q_{st}),
\end{equation}
and such that for all $Q \in \cP(D(\R^+,\cX))$ satisfying
(\ref{eqdiff1}), for all $n$, $\lim_{ t \to +\infty} Q^n (t)=
Q_{\st}^n$.
\end{itemize}
Then the asymptotic independence of the particles  also holds in the
stationary regime:

\begin{theorem}\label{theost}
Under Assumptions A0-A10, for all subsets $I\subset \N$ of finite
cardinal $\vert I\vert$,
\begin{equation*}\label{resEq}
\lim_{N\to \infty} {\cal L}_\st\left((q_i^N(.))_{i\in I}\right)
=Q_\st^{\otimes \vert I\vert} \quad \hbox{weakly in } {\cal
P}(D(\R^+,\cX)^{\vert I\vert}).
\end{equation*}
\end{theorem}

\section{Proof of Theorems \ref{main}, \ref{theoex} and
\ref{theost}}\label{sec:proofs}

We use the following notation extensively:
\begin{equation}
\label{eq:As} A_i ^ {N,s} (k) = \{ \hbox{$s$-transition occurs for
the particle $i$ between $k$ and $k+1$} \}.
\end{equation}
By definition, we have: $$ \PP (A_i ^ {N,s} (k) | \cF_k ) = \frac
1 N F^N_s \left(q_i ^N (\frac  k N) , \mu^N (\frac k N),
Z_i^N(k)\right).$$ We also recall the notation
\begin{equation}\label{eq:Ai}
A^N_i (k) = \{ \hbox{a transition occurs for particle $i$
between times $k$ and $k+1$}\}.
\end{equation}
We have: $A^N_i (k)=\cup_{s\in\cS}A_i^{N,s}(k)$.

\subsection{Proof of Theorems \ref{main} and \ref{theoex}}

By Proposition 2.2. in Sznitman \cite{sznitman}, Theorem \ref{main}
is equivalent to
\begin{equation}
\label{eq:sz}
\lim_{N \to \infty} \cL(\mu^N) = \delta_Q \quad \hbox{
weakly in } {\cal P}({\cal P}(D(\R^+,{\cal X}))).
\end{equation}
To establish (\ref{eq:sz}), we first prove the tightness of the
sequence ${\cal L}(\mu^N,\Gamma^N)$. We then show that any accumulation point
of the previous sequence is the unique solution of a martingale
problem. 
This requires idea from Theorem 2.1 and Example 2.3 in \cite{kurtzaverage}.

\subsubsection{Step 1 : Relative Compactness}

First we check that the sequence ${\cal L}(\mu^N)$ is tight in ${\cal
P}({\cal P}(D(\R^+,\cX)))$. Thanks again to Sznitman \cite{sznitman}
Proposition 2.2, this a consequence of the tightness of ${\cal L}(q_1^N(.))$ in ${\cal P}(D(\R^+,\cX))$; i.e. of A8.
By Prohorov's theorem ${\cal L}(\mu^N)$ is relatively compact. By Lemma 1.3 in \cite{kurtzaverage}, $\Gamma^N$  is relatively compact
because of the compact containment hypothesis A7.  It follows that the
sequence ${\cal L}(\mu^N,\Gamma^N))$ is relatively compact.

\subsubsection{Step 2 : Convergence to the solution of a martingale problem}

We will follow the Step 2 in Graham \cite{graham}. We show that any
accumulation point of ${\cal L}(\mu^N,\Gamma^N)$ satisfies a certain
martingale problem. For $f\in L^{\infty}(\cX)$, the bounded and
forcibly measurable functions of $\cX\to \R$. For each $s \in \cS$,
we define
$$
f^s ( x) = f(s(x)) - f(x).
$$
Now, for $f\in L^{\infty}(\cX)$, $T \geq 0$, and $t\le T$,
\begin{eqnarray}
f(q_i^N(T)) &-& f(q_i^N(t))=\sum_{k = [Nt]}^{[NT]-1}\left(f(q_i^N({(k+1)\over N}-f(q_i^N({k\over N})) \right)\nonumber\\
&=&\sum_{s \in \cS} \sum_{k=[Nt]}^{[NT]-1} f^s(q_i^N ({k\over N}))
\Bigm( \chi\{� A_i ^ {N,s} (k) \} - \PP ( A_i ^ {N,s} (k) | \cF_k)
 \Bigm) \nonumber \\
&& \quad\quad\quad  + \sum_{s \in \cS} \sum_{k=[Nt]}^{[NT]-1} f^s(q_i^N
({k\over N}))\PP ( A_i ^ {N,s} (k) | \cF_k) \label{eq:partiT}.
\end{eqnarray}
Then we define $M^{f,N}_i (t)= \sum_{s \in \cS} M^{f,N,s}_i(t)$ with
\begin{equation}
M^{f,N,s}_i(t)=\sum_{k=0}^{[Nt]-1} f^s(q_i^N ({k\over N})) \left(
\chi\{ A_i ^ {N,s} (k) \}-  \PP ( A_i ^ {N,s} (k) |
\cF_k)\right)\label{mart}
\end{equation}
and
\begin{eqnarray*}
{\cal G}^{N,s}_i f(k)= f^s(q_i^N ({k\over N}))F^{N}_s \left(q_i^N
( {k\over N} ), \mu^N({k\over N}) , Z^N_i(k)\right).
\end{eqnarray*}
So that, we may rewrite Equation (\ref{eq:partiT}) as
\begin{eqnarray}\label{useful}
f(q_i^N(T))&-&f(q_i^N(t))= M^{f,N}_i(T) -M^{f,N}_i(t) + \frac 1 N
\sum_{k=[Nt]}^{[NT]-1}\sum_{s \in \cS} {\cal G}^{N,s}_i f(k)\nonumber\\
&=&M^{f,N}_i(T)- M^{f,N}_i(t)\nonumber\\
&& +  \int_t^T \sum_{s\in \cS}f^s(q_i^N(u))F_s^N(q_i^N(u),\mu^N(u),Z^N_i(u))du.
\end{eqnarray}

The proof of the following lemma is given at the end of this
section.
\begin{lemma}
\label{le:doob}  $M^{f,N}_i(t)$ defined at (\ref{mart}) is a
square-integrable martingale. There exists $C>0$ such that the
Doob-Meyer brackets $\langle M^{f,N}_i, M^{f,N}_i \rangle_t \leq C
t \|f\|^2_{\infty}$ and for $i\neq j$, $| \langle M^{f,N}_i,
M^{f,N}_j \rangle_t | \leq C t \| f \|^2 _{\infty}
\max(\rho_N,1/N)$.
\end{lemma}

Now assume that Lemma \ref{le:doob}  holds, and let $\Pi^{\infty}$ be an accumulation point of $\cL(\mu^N,\Gamma^N)$.
Let $(\mu,\Gamma)$ be a random variable taking values in
$ {\cal P}(D(\R^+,\cX))\times\ell_m(\cZ\times\cX)$ having distribution $\Pi^{\infty}$
which is adapted to a complete filtration $\cF_t$ in the sense that for each $t$,
$\Gamma([0,t],x,z)$ is $\cF_t$-measurable.
By continuity $\Gamma([0,t],x,\cZ)=\int_0^t\mu^x(s)ds$, where $\mu^{x}(s) = \mu(s)(\{x\})$.
By Lemma 1.4 in \cite{kurtzaverage} there exists an $\cF_t$-predictable $\cP(\cX,\cZ)$ valued process $\gamma$
such that $\Pi^{\infty}$-almost surely,
$$\Gamma([0,t],x,z)=\int_{0}^t\gamma_u(x,z)du.$$ 
Define the Radon-Nikodym derivative:
$$\gamma_{(t,x,\Gamma)}(z)=\frac{\Gamma(dt,x,z)}{\Gamma(dt,x,\cZ)}=\frac{\Gamma(dt,x,z)}{\mu^x(t)dt}.$$
Clearly $\gamma_{(t,x,\Gamma)}(z)=\gamma_t(x,z)/\gamma_t(x,\cZ)$  $\Pi^{\infty}$-almost surely.

\begin{lemma}\label{later} We have:
$$\gamma_{(t,x,\Gamma)} = \pi_{\mu(t),x}.$$
\end{lemma}
\bp  Define $\Gamma^N_k (x,z) =\frac{1}{N}\sum_{i=1}^{N}\chi\{X^N_i(k)=x,Z^N_i(k)=z\}  $, we have
\begin{eqnarray}
\lefteqn{\Gamma^N([0,t],x,z)=\frac{1}{N}\sum_{k=0}^{[Nt]}\Gamma^N_k (x,z)   }\nonumber\\
&=&\frac{\Gamma^N_0 (x,z)}{N}  + \frac{1}{N} \hspace{-2pt} \sum_{k=0}^{[Nt]-1}E( \Gamma^N_{k+1}(x,z)|\cF_{k}) \nonumber \\
&&+ \frac{1}{N} \hspace{-2pt} \sum_{k=0}^{[Nt]-1} \hspace{-7pt}\left(\Gamma^N_{k+1} (x,z) - E( \Gamma^N_{k+1}(x,z)|\cF_{k}) \right). \label{eq:gammamart}
\end{eqnarray}
The first term in the above expression goes to $0$ as $N$ goes to infinity. The third term is a mean zero martingale. From Dynkin formula, we have
\begin{eqnarray}
\lefteqn{ E \left( \frac{1}{N} \hspace{-2pt} \sum_{k=0}^{[Nt]-1} \hspace{-7pt}\left(\Gamma^N_{k+1} (x,z) - E( \Gamma^N_{k+1}(x,z)|\cF_{k}) \right) \right) ^2  }\nonumber\\
&=& \frac{1}{N^2} \sum_{k=0}^{[Nt]-1} E  \left(\Gamma^N_{k+1} (x,z) - E( \Gamma^N_{k+1}(x,z)|\cF_{k})  \right) ^2\leq  \frac{ t }{N}  .\nonumber
\end{eqnarray}
 The second term in (\ref{eq:gammamart})
 is equal to
\begin{eqnarray}
\lefteqn{\frac{1}{N}\sum_{k=0}^{[Nt]-1}E( \Gamma^N_{k+1}(x,z)|\cF_{k}) }\nonumber\\
&=& \frac{1}{N}\sum_{i=1}^{N}\frac{1}{N}\sum_{k=0}^{[Nt]-1}E\left(\chi\{X^N_i(k+1)=x,Z^N_i(k+1)=z\}|\cF_{k}\right) \nonumber\\
&=&\frac{1}{N}\sum_{i=1}^{N} \frac{1}{N}\sum_{k=0}^{[Nt]-1}\sum_y \chi\{X^N_i(k)=x ,Z^N_i(k)=y \} K^N_{\nu^N(k),x}(y,z)\cdot (1-\frac{1}{N}) \nonumber \\
& + & \hspace{-6pt} \frac{1}{N}\sum_{i=1}^{N} \frac{1}{N}\sum_{k=0}^{[Nt]-1}\hspace{-7pt} \sum_{w:s(w)=x}  \hspace{-8pt} \chi\{X^N_i(k)=w \} P ( A_i ^{N,s} (k)  ,  Z^N_i(k+1) = z   |\cF_{k} ).  \label{bigone}
\end{eqnarray}
Note that from (\ref{jumps}) $$\hspace{-5pt}\sum_{w:s(w)=x}\hspace{-8pt}  \chi\{X^N_i(k)=w  \}P( A_i ^{N,s} (k)  ,  Z^N_i(k+1) = z   |\cF_{k}) \leq P\left( A_i ^{N} (k) |\cF_{k}\right) \leq \frac 1 N .$$ Thus, as $N\to\infty$ the only important term in (\ref{bigone}) is the first sum and it is equivalent to:
\begin{eqnarray*}
\lefteqn{\frac{1}{N}\sum_{i=1}^{N} \frac{1}{N}\sum_{k=0}^{[Nt]-1} \sum_y \chi\{X^N_i(k)=x ,Z^N_i(k)=y \} K^N_{\nu^N(k),x}(y,z) }\\
&=&\frac{1}{N}\sum_{i=1}^{N} \frac{1}{N} \sum_{k=0}^{[Nt]-1} \sum_y \chi\{q^N_i(\frac k N)=x ,Z^N_i(k)=y \} K^N_{\mu^N(\frac k N),x}(y,z) \\
&=& \int_{0}^{[Nt]/N}\sum_{y\in\cZ}\Gamma^N(du,x,y)K^N_{\mu^N(u),x}(y,z)du\\
&\to& \int_{0}^t\sum_{y\in\cZ}\Gamma(du,x,y)K_{\mu(u),x}(y,z)du
\end{eqnarray*}
as $N\to\infty$ (by Assumptions A3-A4). Therefore, our calculation gives,
$$E\left(\Gamma([0,t],x,z)-\int_{0}^t\sum_{y\in\cZ}\Gamma(du,x,y)K_{\mu(u),x}(y,z)du\right)^2=0.$$
It follows that $\Gamma([0,t],x,z)=\int_{0}^t\sum_{y\in\cZ}\Gamma(du,x,y)K_{\mu(u),x}(y,z)du$ almost surely
and hence that $\gamma_t(x,z)=\sum_{y\in\cZ}\gamma_t(x,y)K_{\mu(t),x}(y,z)$ almost everywhere in $t$  $\Pi^{\infty}$-almost surely.
However for a given $\mu(t)$ and $x$, by Assumption A5  there is a unique solution to the above which is a probability; i.e. for all $z \in \cZ$,
$\gamma_t(x,z)/\gamma_t(x,\cZ)=\pi_{\mu(t),x}(z)$.
\ep

\begin{lemma}
\label{le:G(R)} $\mu$ satisfies a non-linear martingale problem
starting at $Q_0$. Specifically, for all $f\in L^\infty(\cX)$,
\begin{eqnarray}\label{martprob}
M^f (T)=f(X(T))-f(X(0))-\int_0^T {\cal G}f(X(u),\mu(u))du
\end{eqnarray}
is a $\mu$-martingale, where $X=(X(t))_{t\ge 0}$ denotes a canonical
trajectory in $D(\R^+,\cX)$,  $\mu(0)=Q_0$, $\Pi^{\infty}$-a.s. and
$${\cal G}f(x,\mu(t))=\sum_s f^s(x)\overline{F}_s(x,\mu(t)).$$
\end{lemma}

\bp The proof is similar to Step 2 of Theorem 3.4 of Graham
\cite{graham} or of Theorem 4.5 of Graham and M\'el\'eard
\cite{graham97}. However, here our assumptions are weaker so we
detail the proof.

From Lemma 7.1 in Ethier and Kurtz \cite{ethier}, the projection map
$X\mapsto X(t)$ is $\mu$-a.s. continuous for all $t$ except perhaps in
at most a countable subset $D_\mu$ of $\R_+$. Further it is shown
easily that $D=\{t\in \R_+ :\Pi^{\infty}(\{\mu:t\in D_\mu\})>0\}$ is at
most countable (see the argument in the proof of Theorem 4.5 of
Graham and M\'el\'eard \cite{graham97}).

Take $0\leq t_1<t_2<\cdots t_k\leq t <T$ outside $D$ and $g\in
L^{\infty}(\cX^k)$. Take $f\in L^{\infty}(\cX)$. The map $G:{\cal
P}(D(\R^+,\cX))\to \R$ defined by
\begin{eqnarray*}
R\mapsto\langle\left( f(X(T))-f(X(t))-\int_t^T 
{\cal G}f(X(u),\mu(t))du\right)g(X(t_1),\ldots ,X(t_k)),R\rangle
\end{eqnarray*}
is $\Pi^{\infty}$-a.s. continuous. We will prove that
\begin{equation}
\label{eq:G(R)} \Pi^{\infty}\hbox{-a.s, } \quad G(\mu) = 0
.\end{equation} Now assume (\ref{eq:G(R)}) holds for arbitrary
$0\leq t_1<t_2<\cdots t_k\leq t <T$ outside a countable set $D$ and
$g\in C_b(\cX^k)$. It implies that for all $A \subset \mathcal F_t$,
$\langle M^f(T) \ind_A , \mu \rangle = \langle M^f(t) \ind_A,
\mu\rangle$. Therefore, by definition, $M^f(t)$ is a $\mu$-martingale
and $\mu$ satisfies the non-linear martingale problem
(\ref{martprob}).

It remains to prove (\ref{eq:G(R)}).  Let $\Pi^N$ be the law of
$(\mu^N,\Gamma^N)$, we write :
\begin{eqnarray*}
\lefteqn{\langle G , \Pi ^N \rangle
 =   G(\frac 1 N \sum_{i=1} ^ N \delta_{q_i ^N}) }\nonumber \\
& = &  \frac 1 N \sum_{i=1} ^ N  \left( f(q_i ^ N (T)) -
f(q_i ^ N (t))\right)g^{N}_i 
 -  \frac 1 N \sum_{i=1} ^ N  \left(\int_t^T {\cal G}f(q_i ^N (u),\mu(u))du\right)g^{N}_i, \nonumber
\end{eqnarray*}
where $g^{N}_i = g(q_i ^ N  (t_1),\ldots ,q_i ^ N  (t_k))$.
From (\ref{useful}),
\begin{eqnarray*}
\lefteqn{\langle G , \Pi ^N \rangle
 =  \frac 1 N \sum_{i=1} ^ N \left(M^{f,N}_i(T)-M^{f,N}_i(t)\right)g^{N}_i}\\
&& + \frac 1 N \sum_{i=1} ^ N \left(\int_t^T \sum_{s\in \cS}f^s(q_i^N(u))F_s^N(q_i^N(u),\mu^N(u),Z^N_i(u))du\right.\\
&& \left.- \int_t^T {\cal G}f(q_i ^N (u),\mu(u))du\right)g^{N}_i.
\end{eqnarray*}
Hence,
\begin{eqnarray*}
\lefteqn{\EE |\langle G , \Pi ^N \rangle|\leq  \EE  \left| \frac 1 N \sum_{i=1} ^ N (M^{f,N}_i (T) -
M^{f,N}_i (t) ) g^N_i  \right|}\nonumber\\
& &+ \EE  \left|\frac 1 N \sum_{i=1} ^ N \left(\int_t^T \sum_{s\in \cS}f^s(q_i^N(u))F_s^N(q_i^N(u),\mu^N(u),Z^N_i(u))du\right.\right.\\
&&\left.\left. - \int_t^T {\cal G}f(q_i ^N (u),\mu(u))du\right)g^N_i  \right|\nonumber\\
& \leq & \hbox{I} + \hbox{II} \label{eq:G2},
\end{eqnarray*}

Using exchangeability and the Cauchy-Schwartz inequality, we obtain:
\begin{eqnarray*}
\hbox{I}^2 & \leq &   \frac {\|g \|^2_{\infty} }{N} \EE
\left(M^{f,N}_1 (T) - M^{f,N}_1 (t)\right)^2 \\ && \quad  +
\frac{N-1}{N} \EE\Bigm( (M^{f,N}_1 (T) - M^{f,N}_1
(t))g^{N}_1(M^{f,N}_2 (T) - M^{f,N}_2 (t))g^{N}_2\Bigm).
\end{eqnarray*}
Lemma \ref{le:doob} implies that $\hbox{I}$ tends to $0$. 

Next, $\hbox{II}^2$ is less than or equal to
\begin{eqnarray*}
\lefteqn{ \|g \|^2_{\infty} \EE \Bigg(\frac 1 N \sum_{i=1} ^ N \bigg(\int_t^T \sum_{s\in \cS}f^s(q_i^N(u))F_s^N(q_i^N(u),\mu^N(u),Z^N_i(u))du }\\
& &\quad\quad - \int_t^T {\cal G}f(q_i ^N (u),\mu^N(u))du\bigg)\Bigg)^2\\
&\leq&\|g \|^2_{\infty}\EE \Bigg(\int_t^T \sum_{s\in \cS}\sum_{x,z}f^s(x)F_s^N(x,\mu^N(u),z)\Gamma^N(du,x,z)\\
& &\quad\quad - \int_t^T \sum_{x}{\cal G}f(x,\mu^N(u))\Gamma^N(du,x,\cZ)\Bigg)^2.
\end{eqnarray*}
However, as $N\to\infty$,
\begin{eqnarray*}
\lefteqn{\int_t^T \sum_{s\in \cS}\sum_{x,z}f^s(x)F_s^N(x,\mu^N(u),z)\Gamma^N(du,x,z)}\\
&\to&\int_t^T \sum_{s\in \cS}\sum_{x,z}f^s(x)F_s(x,\mu(u),z)\Gamma(du,x,z),
\end{eqnarray*}
and using Lemma \ref{later}
\begin{eqnarray*}
\lefteqn{ \int_t^T\sum_{x} {\cal G}f(x, \mu ^N (u) ) \Gamma^N(du,x,\cZ) }\\
&=&\int_t^T \sum_s\sum_{x }f^s(x)F_s(x,\mu^N(u),z)\gamma_{(t,x,\Gamma)}(z)\Gamma^N(du,x,\cZ)\\
&\to&\int_t^T \sum_s\sum_{x,z}f^s(x)F_s(x,\mu(u),z)\gamma_{(t,x,\Gamma)}(z)\Gamma(du,x,\cZ)\\
&=&\int_t^T \sum_{s\in \cS}\sum_{x,z}f^s(x)F_s(x,\mu(u),z)\Gamma(du,x,z).
\end{eqnarray*}
Consequently  $\hbox{II}^2 \to 0$ as $N\to\infty$.

 Hence,
from (\ref{eq:G2}) and Fatou's Lemma, $\langle |G| , \Pi ^\infty
\rangle \leq \lim_N \langle |G| , \Pi ^N \rangle = 0$ and thus
$\Pi^{\infty}$-a.s, $G(\mu) = 0$, (\ref{eq:G(R)}) is proved.

To conclude the proof of Lemma \ref{le:G(R)}, note that the
continuity of $X\to X(0)$ implies $\mu(0)=Q_0$,
$\Pi^{\infty}$-a.s..\ep

\subsubsection{Step 3 : Uniqueness of the solution of martingale problem}
We now show the solution to (\ref{martprob}) is unique. Here, we
will use Proposition 2.3 in Graham \cite{graham} (which is an
extension of Lemma 2.3 in Shiga and Tanaka \cite{shiga}) to show
uniqueness. We remark that 
${\cal G}f(x,\alpha)=\int_{\cX}(f(y)-f(x))J_{x,\alpha}(dy)$ where
$$J_{x,\alpha} = \sum_{s \in \cS} \overline F_s (x,\alpha) \delta_{s(x)}.$$
Next, $\|J_{x,\alpha}\| =  \sum_{s \in \cS} \overline F_s (x,\alpha)
= 1$ and  $\|J_{x,\alpha}-J_{x,\beta}\| = \sup |\int_{\cX} \varphi(y
) J_{x,\alpha}(dy) - \int_{\cX} \varphi(y ) J_{x,\beta}(dy)|$ where
the supremum is over the functions $\varphi \in L^{\infty} (\cX)$
with $\|\varphi\|_{\infty} \leq 1$.
\begin{eqnarray*}
 |J_{x,\alpha}(\varphi) &-& J_{x,\beta} (\varphi)| =  \Bigm|\sum_{s \in \cS} \varphi(s(x)) \Bigm(\overline F_s
(x,\alpha) - \overline F_s
(x,\beta)\Bigm) \Bigm| \\
 & = & \Bigm| \sum_{s \in \cS} \varphi(s(x)) \Bigm( \int_{\cZ}  F_s
(x,\alpha,z) \pi_{\alpha,x} (dz)  - \int_{\cZ} F_s (x,\beta,z)
\pi_{\beta,x} (dz)\Bigm) \Bigm| \\
 &\leq& \Bigm| \sum_{s \in \cS} \varphi(s(x)) \Bigm(\int_{\cZ}  F_s
(x,\alpha,z) \pi_{\alpha,x} (dz)  - \int_{\cZ} F_s (x,\alpha,z)
\pi_{\beta,x} (dz) \Bigm) \Bigm|  \\
&&+ \Bigm| \sum_{s \in \cS} \varphi(s(x)) \Bigm(\int_{\cZ} F_s
(x,\alpha,z)
\pi_{\beta,x} (dz) - \int_{\cZ} F_s (x,\beta,z) \pi_{\beta,x} (dz)\Bigm) \Bigm|\\
& \leq & \hbox{I} + \hbox{II}.
\end{eqnarray*}
By Fubini's Theorem, \begin{eqnarray*} \hbox{I} & = &
\Bigm|\int_{\cZ} \sum_{s \in \cS} \varphi(s(x))   F_s (x,\alpha,z)
\pi_{\alpha,x} (dz) - \int_{\cZ} \sum_{s \in \cS} \varphi(s(x))
F_s (x,\alpha,z)
\pi_{\beta,x} (dz)\Bigm| \\
& \leq & \|\sum_{s \in \cS}  \varphi(s(x)) F_s (x,\alpha,\cdot)
\|_{\infty} \| \pi_{\alpha,x} - \pi_{\beta,x} \|.
\end{eqnarray*}
Since $ |\varphi(s(x))| \leq 1$,  $F_s (x,\alpha,z) \geq 0$ and
$\sum_{s \in \cS} F_s (x,\alpha,z) = 1,$
$$
\|\sum_{s \in \cS} \varphi(s(x)) F_s (x,\alpha,\cdot) \|_{\infty}
\leq 1.
$$
Thus applying Assumptions A4-A6, we deduce:
$$
\hbox{I} \leq  \| \pi_{\alpha,x} - \pi_{\beta,x} \| \leq C
\|\alpha - \beta\|.
$$
Using Assumption A2,
\begin{eqnarray*}
\hbox{II} & \leq &  \int_{\cZ} \sum_{s \in \cS}  | F_s (x,\alpha,z)
 -  F_s (x,\beta,z)  |\pi_{\beta,x} (dz)\\
 & \leq & C \|\alpha - \beta\|.
\end{eqnarray*}
So finally, we have checked that:
$$
\| J_{x,\alpha} - J_{x,\beta} \| \leq C \|\alpha - \beta\|.
$$
We then use Proposition 2.3 in Graham \cite{graham} to establish the
solution to the martingale problem (\ref{martprob}) is unique.

\subsubsection{Step 4 : Weak convergence and Evolution equation}

In the three first steps we have proved that $\cL(\mu^N)$ converges
weakly to $\mu=\delta_Q$, where $Q$ is the unique solution of the
martingale problem (\ref{martprob}) starting at $Q_0$. 

We can now identify the evolution equation satisfied by $Q$.  Since
$Q$ satisfies the martingale problem then $(Q(t))_{t\geq 0}$ solves
the non-linear Kolmogorov equation derived by taking the
expectations in (\ref{martprob}):
\begin{eqnarray}\label{crochet}
\langle f,Q(T) \rangle-\langle f,Q(0)\rangle=\int_0^T\langle {\cal
G}f(\cdot ,Q(t)),Q(t)\rangle dt.
\end{eqnarray}

Applying (\ref{crochet}) to $f=\ind_{x_n}$ for all $n$, we get the
set of differential equations (\ref{eqdiff1}). It immediately follows that $\Gamma$ is also deterministic and $\Gamma(dt,x,z)=dt\cdot Q^x(t)\cdot \pi_{Q(t),x}(z)$ almost surely.

\subsubsection{Proof of Lemma \ref{le:doob}}

First, $M^{f,N}_i(t)$ is a square-integrable martingale by the
Dynkin formula. Recall that $A^{N,s}_i (k)$ is defined in Equation
(\ref{eq:As}) and that $A^N_i (k) = \cup_{s \in \cS} A^{N,s}_i
(k)$. In the sequel, $\EE_{\cF_k} [.]$ will denote $\EE[. |
\cF_k]$. With this notation, $\EE_{\cF_k} [\ind_{A^N_{i}(k)}] =
1/N$, and we can rewrite Equation (\ref{mart}) as:
\begin{equation*}
M^{f,N}_i(t)=\sum_{k=0}^{[Nt]-1}\sum_{s \in \cS} f^s(q_i^N({k\over
N}))\Bigm( \chi\{A^{N,s}_i (k)\}  - \EE_{\cF_k}
\chi\{A^{N,s}_i (k)\} \Bigm).\label{mart2}
\end{equation*}
To prove Lemma \ref{le:doob}, we first need to compute $\EE [M^{f,N}_1(t)
M^ {f,N}_2(t)]$. Since $(M^{f,N}_i(t))_{t \in \R^+}$ is a
martingale this product is equal to:
\begin{eqnarray*}
\lefteqn{ \EE [M^{f,N}_1(t) M^ {f,N}_2(t)] }\\
& = & \sum_{k = 0} ^ {[Nt]-1}
\sum_{s,s' \in \cS} \EE f^s(q_1^N({k\over N}))\Bigm(  \chi\{A^{N,s}_1 (k)\}  - \EE_{\cF_k}
\chi\{A^{N,s}_1 (k)\} \Bigm) \\
&& \quad  \times  f^{s'}(q_2^N({k\over N}))\Bigm( \chi\{A^{N,s'}_2
(k)\} - \EE_{\cF_k} \chi\{A^{N,s'}_{2}(k)\} \Bigm).\end{eqnarray*}
 Now, let
\begin{eqnarray*}
I^N_k &=&  \sum_{s,s' \in \cS} \EE \left[ f^s(q_1^N({k\over
N}))(  \chi\{A^{N,s}_1 (k)\}  - \EE_{\cF_k}
\chi\{A^{N,s}_1 (k)\} )\right.\\
&&\quad\quad\quad \quad\quad\times \left. f^{s'}(q_2^N({k\over N}))(
 \chi\{A^{N,s'}_2
(k)\} - \EE_{\cF_k} \chi\{A^{N,s'}_{2}(k)\} )\right] \\
& = & \sum_{s,s' \in \cS} \EE \left[ f^s(q_1^N({k\over N}))  f^{s'}(q_2^N({k\over N})) \right.\\
&&\times \left.\Bigm( \EE_{\cF_k} [ \chi\{A^{N,s}_1 (k)\} \chi\{A^{N,s'}_2 (k)\}] - \EE_{\cF_k} [
\chi\{A^{N,s}_{1}(k)\} ]\EE_{\cF_k} [
\chi\{A^{N,s'}_{2}(k)\}]  \Bigm)\right].
\end{eqnarray*}
Notice that
\begin{eqnarray*}
\lefteqn{ \Bigm| \sum_{s,s' \in \cS} f^s(q_1^N({k\over N}))  f^{s'}(q_2^N({k\over N}))  \EE_{\cF_k} [ \chi\{A^{N,s}_1 (k)\} \chi\{A^{N,s'}_2 (k)\} ] \Bigm| }\\
& &\leq  4 \| f \|^2_{\infty} \EE_{\cF_k} [ \sum_{s,s'}  \chi\{A^{N,s}_1 (k)\} \chi\{A^{N,s'}_2 (k)\} ]    \\
& &\leq  4 \| f \|^2_{\infty}  \PP (  A^{N}_1 (k) A^{N}_2 (k) | \cF_k).
\end{eqnarray*}
Analogously, we also have: 
\begin{eqnarray*}
\lefteqn{ \Bigm| \sum_{s,s' \in \cS}
f^s(q_1^N({k\over N}))  f^{s'}(q_2^N({k\over N})) \EE_{\cF_k} [
\chi\{A^{N,s}_{1}(k)\} ]\EE_{\cF_k} [ \chi\{A^{N,s'}_{2}(k)\}]
\Bigm| }\\ 
& &\leq  4 \| f \|^2_{\infty}  \PP (A^{N}_1 (k)|
\cF_k)\PP(A^{N}_2 (k)| \cF_k).
\end{eqnarray*}
Therefore from
(\ref{cancellation}),  $ | I^N_ k  | \leq  8 \| f\|_{\infty} ^ 2
\max(\rho_N/N,1/N^2)$ and
$$ | \EE [M^{f,N}_1(t) M^ {f,N}_2(t)] | \leq 8 \| f\|_{\infty} ^ 2 t \max(\rho_N,1/N) .$$ Similarly, we obtain
\begin{eqnarray*}
\EE \left[\left(M^{f,N}_1(t) \right)^2 \right] & = & \sum_{k = 0} ^ {[Nt]-1}
 \EE\left( \sum_{s \in \cS}  f^s(q_1^N({k\over N}))\Bigm( \chi\{A^{N,s}_1
(k)\}  - \EE_{\cF_k} \chi\{A^{N,s}_{1}(k)\} \Bigm) \right)^2 \\
& \leq & \sum_{k = 0} ^ {[Nt]-1} 8 \| f \|^2_{\infty} \PP ( A^{N}_1 (k) ) \\
& \leq &  8 \| f \|^2_{\infty} t,
\end{eqnarray*}
and the lemma follows.
\ep

\subsection{Proof of Theorem \ref{theost}}

Assume that $((q_i ^N (0))_i,Z^N)$ represents the system of $N$
particles in stationary regime. Then by symmetry, $(q_i ^N (0))_i$
is exchangeable. Define $\Pi^N = \cL(\mu^N,\Gamma^N)$. We cannot apply directly
Theorem \ref{main} since we do not know whether a converging subsequence of
$\mu^{N} (0)$ converges weakly toward a deterministic limit.

We now circumvent this difficulty. By Assumption A9, as in Step 1 in the proof of
Theorem \ref{main}, we deduce from Sznitman \cite{sznitman}
Proposition 2.2, that $\mu^N$ is tight in  $ {\cal P}(D(\R^+,\cX)) $ and $\Pi^N$ is tight in  $ {\cal P}(D(\R^+,\cX))\times\ell_m(\cZ\times\cX)  $.
 Let $Q$ in ${\cal P}(D(\R^+,\cX))$ be in the
support of $\Pi^{\infty}=(\mu^{\infty},\Gamma^{\infty})$, an accumulation point of $\Pi^N$. We can
prove similarly that Lemma \ref{le:G(R)} still holds for $Q$.

By Step 3 of Theorem \ref{main}, the solution of the martingale
problem is unique and $Q$ solves it with initial condition $Q(0)$.
The stationarity implies that $\mu^N (t)$ and $\mu^N (0)$ are equal.
Note also that outside a countable set $D$, the mapping $X \mapsto
X(t)$ is continuous. So if $ t\notin D$,
$Q(t)=\mu^{\infty}(t)=\mu^{\infty} (0)=Q(0)$.  However, by Assumption A9,
$\lim_{t \to + \infty} Q(t) = Q_{st}$. Therefore $\mu^\infty (0) =
\delta_{Q_{st}}$ and $Q(0)=Q_{st}$.

Theorem \ref{theost} is then a consequence of Theorem \ref{main}.

\section{A uniform domination criterion} \label{sec:HYP}

In this section we discuss the Assumptions A0-A9 made on the
particle system. Assumptions A0-A6 are natural and can be checked
directly.  The additional assumptions
A9 and A10 needed to derive the mean field limit in the stationary
regime may be difficult to check: A9 is a tightness assumption on
the stationary measures and A10 is the global stability of a
differential equation.

In this section we present a new set of assumptions, based on
uniform domination of the transition kernel of the background
process, that is provably sufficient to ensure that Assumptions
A7. The new assumptions are defined as follows:

\begin{itemize}
\item[A11] There exists a
transition kernel $K$ on $\cZ$ which dominates the kernels
$K^N_{\alpha,x}$. Specifically, let $\preceq$ be a partial order
on $\cZ$ such that $\cK_z = \{�w \in \cZ : w \preceq z \}$ is finite for all $ z \in \cZ$. There exists $K$ such that
for all $N$, $z$, $x$, $\alpha$,
$$
K_{\alpha,x}^N  (z,\cdot) \preceq_{st} K(z,\cdot),
$$
where $\preceq_{st}$ is the stochastic order relation: $P
\preceq_{st} P'$ if for all $z_1 \in \cZ$: $ \sum_{z \succeq z_1}
P(z) \leq \sum_{z \succeq z_1} P'(z)$.
\item[A12] The Markov chain $Z(t)$  with transition kernel $K$ is positive recurrent.
\end{itemize}

\begin{lemma}
Under Assumptions A8 and  A11-A12, A7 holds. 
\end{lemma}

\bp  Because the chain $Z$ is  positive recurrent, the long run proportion of time the chain $Z$ spends
outside a compact set $\cK_z$ is of probability at most $\epsilon/2$ for some $z \in \cZ$, so
\begin{eqnarray*}
\frac{1}{t}E\Gamma^N([0,t] \times \cX \times \cK_z)&=&\frac{1}{Nt}\sum_{k=0}^{[Nt]}(\frac{1}{N}\sum_{i=1}^NE\chi\{Z_i^N(k)\in \cK_z\})\\
&\geq&\frac{1}{Nt}\sum_{k=0}^{[Nt]}P(Z(k)\in \cK_z)\\
&\to &(1-\epsilon/2)
\end{eqnarray*}
as $N \to\infty$
Hence $\liminf_N E\Gamma^N([0,t]\times\cX\times \cK_z)\geq (1-\epsilon/2)t$.

By A8 we know $\mu^N$ is relatively compact and hence tight. By (2.5) in \cite{sznitman}
the tightness of $\mu^N$ is equivalent to the tightness of their intensity measures $I(\mu^N)$ in $\cP(D(\R^+,\cX))$ defined by
$I(\mu^N)(F)=E\mu^N(F)$ for $F\in \cB((D(\R^+,\cX))$, the Borel $\sigma$-algebra associated to the Skorohod topology. Hence for every $\epsilon>0$  there exists a compact set $K_{\epsilon}$
in $D(\R^+,\cX)$ such that
$\inf_N E\mu^N(K_{\epsilon})  \geq 1-\epsilon/2.$
However by Remark 6.4 on page 124 in \cite{ethier}, for each $T>0$ there exists a compact set
$\cK_{x}\subseteq \cX$ such that for all $t\in [0,T]$, $\{�x ( t)  : x (\cdot) \in K_{\epsilon}\} \subseteq \cK_x$.
Hence, $ \mu^N(t)(\cK_x)\geq  \mu^N(K_{\epsilon})$ for all $t\in [0,T]$ and for all $N$. Consequently
$\inf_NE\mu^N(t)(\cK_{x})   \geq 1-\epsilon/2$ for all $t\in [0,T]$.
However, for each $t$,
\begin{eqnarray*}
\frac{1}{t}E\Gamma^N([0,t],\cK_x,\cZ)&=&\frac{1}{t}\int_0^{t}E\mu^N(u)(\cK_x)du \geq 1-\epsilon/2.
\end{eqnarray*}
by the above. We conclude A7 holds with $\cK = \cK_x \times \cK_z$.

 \ep

%Part of Assumption A8 may be verified using the next lemma.
%\begin{lemma}
%\label{le:tightN} Under Assumptions A10-A13, if in addition,
%$(\mu^N,Z^N)$ is an ergodic Markov chain for all $N$ and if for all
%$z_1 \in \cZ$, $\{ z : z \preceq z_1\}$ is finite then $\cL_{st}
%(Z^N)$ is tight.
%\end{lemma}

%\bp As in the proof of Lemma \ref{le:renewal}, if $Z$ is a Markov
%chain with kernel $K$, $ \PP ( Z^N (k) \succeq t | Z^N (0) = z_1 )
%\leq \PP ( Z(k) \succeq t | Z(0) = z_1 )$. Letting $k$ tends to
%infinity, using the ergodicity of the chains $(\mu^N,Z^N)$ and $Z$,
%we obtain the tightness of $\cL_{st} (Z^N)$. \ep

\section{Application to random multi-access
protocols}\label{sec:example}

We now apply the previous analytical results to study the performance of communication networks where $N$ users share a common resource in a distributed manner. We consider for example Local Area Networks (LANs) which are computer networks with relatively small geographic coverage (an office, a house, a part of a campus), and which constitutes the first crucial component of the Internet. Transmissions in LANs are handled either on a cable (wired LANs) or on a radio channel (wireless LANs, also commonly called WiFi). Here we will focus on wireless LANs (our analysis can be carried out similarly in the case of wired LANs). In wireless LANs, users that are close to each other or that wish to transmit to the same receivers {\it interfere} in the sense that they cannot simultanerously transmit packets succesfully. Two interfering users transmitting simultaneously are said to experience a {\it collision}. A collision is detected by a user at the end of the packet transmission when the corresponding receiver does not acknowledge a successful reception. One of the most challenging problem in computer networking has been to design mechanisms so that interfering users could efficiently and fairly share the resource in a distributed manner. Currently, users willing to transmit packets through a wireless LAN, implement two standardized mechanisms, Carrier Sense Multiple Access (CSMA) and a random back-off algorithm referred to as the Decentralized Coordination Function (DCF), see \cite{ieee}. In this section we aim at analyzing the performance of a general class of mechanisms, including the current CSMA - DCF couple, and at understanding whether current mechanisms perform well or if they still require important improvements.     

In the next subsection, we provide a short description of CSMA and of a class of random back-off algorithms, but also introduce a simple model for interference, and explain why the performance in wireless LANs is difficult to study. In the subsequent subsections, we explain how  the results derived earlier in the paper for particle systems allow us to circumvent this difficulty and explicitly characterize the performance in these networks.

\subsection{Distributed mechanisms and performance in wireless LANs}

\subsubsection{Carrier Sensing mechanisms}

A first mechanism to separate transmissions of interfering users in time is CSMA. Before transmission, each user senses the channel, and should it be busy, it abstains from transmitting. This sensing mechanism may be too simple to capture the actual interference structure of the network (since for example, the sensing is made at the transmitters, whereas interference is experienced at the receivers). Collisions may occur due to {\it hidden} terminals, and a loss of efficiency can be due to {\it exposed} terminals, see e.g. \cite{jiang}. Hidden terminals refer to users whose transmissions interfer at the receiver, but are not able to detect (sense) each other. On the contrary, exposed terminals are users that do not interfere at the receiver, but cannot  simultaneously transmit because they sense each other's transmissions. In this paper, for simplicity, we restrict our attention to a perfect Carrier sensing mechanism, where users sensing each other actually interfere at the receiver (we believe the analysis could be extended with hidden and exposed terminals).

\subsubsection{Random back-off algorithms}   

\begin{figure}[h]\label{fig:dcf}
\centering
\includegraphics[width=8cm]{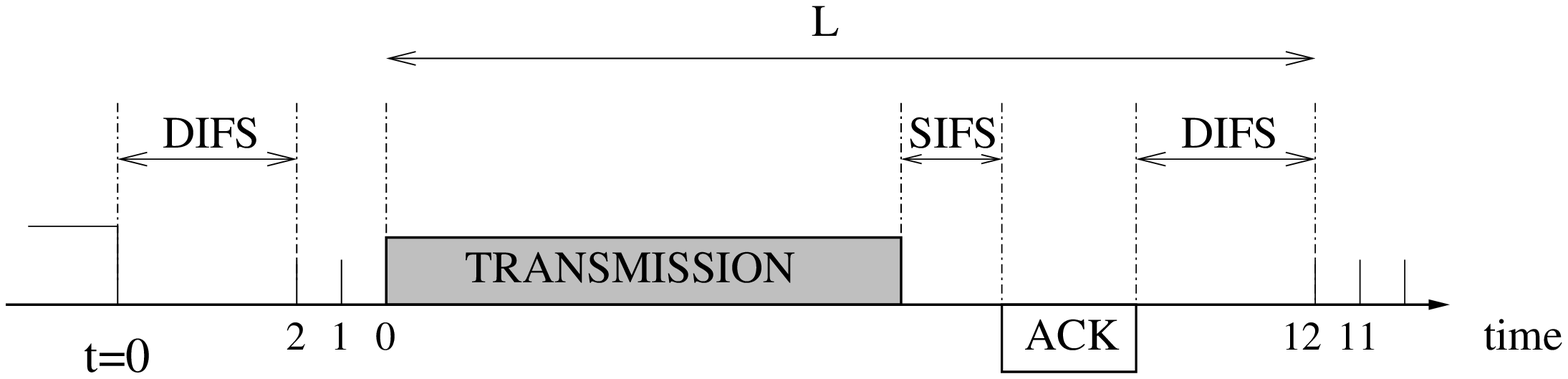}
\caption{User behavior - the case of a successful transmission. Before $t=0$, the channel is sensed busy. At time DIFS, (DCF Inter Frame Space), the user starts decrementing its back-off counter again by one per slot, and transmits when the latter reaches 0. After transmission, the receiver waits for a duration of length SIFS (Short Inter Frame Space) and then sends the packet acknowledgment. After receiving this acknowledgment, the user picks a new back-off counter ($12$ in this case) and waits DIFS before starting decrementing it. Note that the inter-frame spaces are introduced to handle the acknowledgment procedure, and that DIFS $>$ SIFS.}
\end{figure}

Even under a perfect carrier sensing mechanism, collisions cannot be completely avoided if two users start transmitting simultaneously. To further reduce collisions, each user runs (independently of other users) a random back-off algorithm. After each successful transmission or each collision, the user randomly picks a value for its back-off counter according to some distribution on $\mathbb{N}$. This value represents      the number of {\it slots} the channel has to be observed idle before that the user may start transmitting (basically the user decrements its counter by one after sensing the channel idle during one slot). Note that slots have a fixed duration  that does not depend on the user (between 9 and $20$ microseconds in IEEE802.11 standards \cite{ieee}). The details of this mechanism works is exemplified in Figure \ref{fig:dcf}.

A random back-off algorithm specifies how the distribution (or just its mean) of the back-off counter is modified after either a successful transmission or a collision. Currently  the DCF
is a version of the classical binary exponential back-off algorithm: after each successful transmission, a user picks a back-off counter uniformly in $\{0,\ldots, CW_{\min}\}$, and after $m$ successive collisions uniformly in $\{0,\ldots, 2^{m}CW_{\min}\}$\footnote{Note that in the DCF, $m$ is upper bounded by 7.}.

In the following, we assume that the back-off distribution is always geometric (so as to keep a simple Markovian setting), although we could easily generalize the analysis to uniform distributions. With this assumption, each user transmits with a given probability $p$ at the beginning of each idle slot. We consider the following generic way of adapting this probability: first the probability belongs to a countable set ${\cal B}$, after a successful transmission $p$ is updated to $S(p)$, and after a collision $p$ is updated to $C(p)$, where $S(\cdot)$ (resp. $C(\cdot)$) is a decreasing (resp. increasing) mapping from ${\cal B}\to{\cal B}$. We denote by $p_0=\max\{p\in {\cal B}\}$. Finally, we denote by $L$ (in slots) the average duration of a successful packet transmission (including its acknowledgment, see Figure \ref{fig:dcf}), and assume that collisions have average durations equal to $L_c$ that might be different than $L$. Again to keep the formalism simple, we assume that the durations of successful transmissions and collisions are geometrically distributed (a multiple of slots), which again does not constitute a crucial assumption.

\subsubsection{Interference model and user class}

We consider a simple model for interference as follows. First, the $N$ users are classified according to their interference properties, i.e., two
users belong to the same class if they interfere  with (resp. are
interfered with and by) the same set of users. Two users are of the same class if the corresponding links are located in the same geographic region (see for example the network of Figure \ref{firstpic}).
Denote by $\cC$ the set of user classes, and by $\mu_c$ the proportion of users of class $c$. $i\in c$ denotes the fact that user $i$ is of class $c$. Then interference between users of different classes is characterized by the incidence matrix $A$ such that $A_{cd}=1$ if class-$c$ users interfere class-$d$ users, and $A_{cd}=0$ otherwise. 
Note that $A$ is not necessarily symmetric (in the network of Figure \ref{firstpic}, it is symmetric). We denote by $\cV_c=\{d\in\cC: A_{cd}=1\}$ the set of classes of links interfering with class-$c$ links.

We say that the network has full interference if $A_{cd}=1$ for all $c,d$ and has partial interference otherwise.

\subsubsection{Performance metrics}

The performance metrics we aim at analyzing is the long-term throughput (the number of packets successfully transmitted per time unit) achieved by the users of various classes. We denote by $\gamma_c$ the throughput of class-$c$ users.

Deriving expressions for this performance metrics is notoriously difficult. This is due to the inherent interactions between users through interference. A popular approach to circumvent this difficulty consists in decoupling the users, i.e., assuming that the (re)-transmission processes of the various users are mutually independent. This heuristic has been used by Bianchi \cite{bianchi} to capture the performance of wireless LANs with full interference. In this work, we formally justify this approach, and extend it to networks with partial interference. To do so, we apply the mean field analysis derived in
the first part of the paper. In case of full interference, the
network can be modeled as a simple system of particles with no
randomly varying environment (as already noticed in a preliminary work \cite{allerton}). However, to analyze a network with partial interference, the introduction of this varying environment is necessary. As it turns out, the spatial heterogeneity in networks with partial interference may lead to important fairness issues, as mentioned in introduction, and our analysis  explicitly quantifies these issues.

\subsection{Model analysis}\label{subsub:model}

We consider a network of $N$ users as described in the previous subsection. We analyze the system at the beginning of
each slot. Denote by $p_i^N(k)/N$ the probability user $i$ becomes
active at the end of the $k$-th slot, if idle (note that we already
renormalized this probability by $1/N$ to be able to conduct the
asymptotic analysis when $N$ grows large). For all $i,k,N$, $p_i^N(k)\in \cB$.

To capture the network dynamics, we define a process $Z^N=\{Z^N(k), k\ge 0\}$ representing the {\it state} of classes during slot $k$.
 $Z^N_c(k)\in \{0,1,2\}$, where $Z_c(k)=0$ if and only if there is no transmitting user of class $c$, $Z^N_c(k)=1$ if and only if there is one
successfully transmitting user of class $c$ and $Z^N_c(k)=2$ if and
only if there is at least one user of class $c$ currently in
collision with another user in $\cV_c$. Let $\cZ=\{0,1,2\}^{\vert
\cC\vert }$ denote the state space of $Z$. We introduce the {\it
clear-to-send} functions $C_c$ as follows. If $Z^N(k)=z$, a class-$c$
link is clear to send at the end of slot $k$ and $C_c(z)=1$  if
$z_d=0$ for $d\in \cV_c$, otherwise $C_c(z)=0$. 

We show how to model the network as a set of interacting
particles as described in Section \ref{markov}.
\begin{itemize}
\item The particles: the $i$-th user corresponds to the $i$-th particle
with state describing the class of the user and the transmission
probability at the end of the next idle slot
$X_i^N(k)=(c_i,p_i^N(k)) \in \cX=\cC\times\cB$.
\item The environment process: the process $Z^N$ introduced above 
is a simplified version of  an environment process as described in Section \ref{markov}. The evolution of the environment  is determined by the states of all the particules through $\nu^N$. The evolution of the $i$-th particle depends on whether or not the corresponding user senses the channel idle or not, i.e.  $Z^N_i(k)=Z_c^N(k)$ for $i\in c$. 
\end{itemize}

\paragraph{Particle transitions} We first compute the transition
probabilities for the various particles. The set $\cS$ of possible
transitions is composed by two functions, the first one representing
a successful transmission $p\mapsto S(p)$ and the other one
collisions $p\mapsto C(p)$. Note that the class of a particle / user
does not change. Let $\nu^N_c(k)=\frac 1 N \sum_{i=1} ^ N
\delta_{p^N_i(k)}\ind_{c(i)=c}$ and $\nu^N(k)=(\nu^N_c(k))_{c\in
\cal C}$.

Assume that at some slot $k$, the system is in state
$$
((c_i^N(k),p_i^N(k))_{i=1,\ldots
N},\nu^N(k),Z^N(k)))=((c_i,p_i)_{i=1,\dots ,N},\alpha,z).
$$
A class-$c$ user $i$ may have a transition at the end of slot $k$
only if $C_c(z)=1$. In this case it can either initiate a successful
transmission or experience a collision. If $C_c(z)=1$, the event
that none of the users in $c$ transmits at the end of slot $k$ is
given by $D^N_c=\prod_{i\in c}1_{(N U_i>p_i)}$, where the
$U_i$'s are i.i.d. r.v. uniformly distributed on $[0,1]$. The event
that user $i\in c$ accesses  the channel with success at the end
of slot $k$ is given by the indicator:
\begin{equation*}
\chi\{N U_i\le p_i\}C_c(z)\prod_{j\in c, j\neq i}\chi\{N
U_j>p_j\}  \prod_{d\in \cV_{c},d\neq c} \Bigm( C_d(z)
D^N_d + (1-C_d(z)) \Bigm).
\end{equation*}
Averaging the above quantity gives the transition probability
$F_S^N((c,p_i),\alpha,z)/N$ corresponding to a successful
transmission. For all $\alpha\in \cP(\cB)$ and all $f$ $\cB$-valued
functions, define $\langle f,\alpha\rangle=\sum_pf(p)\alpha(p)$.
Moreover let $\alpha_c$ denote the restriction of $\alpha$ to users of class
$c$. Let $I$ denote the identity function. One can readily see
that we have:
\begin{equation}
\label{eq:defFNS} F_S^N ((c,p_i),\alpha,z) = \frac{p_i}{1- p_i/N}
C_c (z)\prod_{ d\in \cV_{c}} \Bigm( C_d(z) ( e^{ \langle
N\log(1-\frac{I}{N}),\alpha_d\rangle} -1 )+1\Bigm).
\end{equation}

Similarly, the event that user $i\in c$ experiences a collision
at the end of slot $k$ is given by the indicator:
\begin{equation*}
\chi\{N U_i\le p_i\}C_c(z)\left(1 -\prod_{j\in c, j\neq i}\chi\{N
U_j>p_j\} \prod_{d\in \cV_{c},d\neq c} \Bigm( C_d(z)
D^N_d + (1-C_d(z)) \Bigm)\right),
\end{equation*}
and the transition probability $F_C^N((c,p_i),\alpha,z)/N$
corresponding to a collision reads:
\begin{equation}
\label{eq:defFNC} F_C^N ((c,p_i),\alpha,z) = p_iC_c (z)\left(1 -
{1\over 1- p_i/N}\prod_{ d\in \cV_{c}}  \Bigm(  C_d(z)( e^{ \langle
N\log(1-\frac{I}{N}),\alpha_d\rangle} -1 )+1\Bigm)\right).
\end{equation}

In order to fit into the scheme to the particle system of Section
\ref{markov}, we need to introduce a virtual transition from $(c,p)$
to $(c,p)$ with transition rate $F_{\emptyset} ^N ((c,p_i),\alpha,z)
= 1 - p_i C_c (z)$. With this virtual transition the sum of the
transition rates sums to $1$. Note that Assumption 0 is satisfied. Since $N \log (1 - x / N) $ converges
to $-x$, we obtain the following expressions for the asymptotic
transition rates, $F_{\emptyset}  ((c,p_i),\alpha,z) = 1 - p_i  C_c
(z)$,
\begin{equation}
\label{eq:defFS} F_S ((c,p_i),\alpha,z) = p_i
  C_c (z) \prod_{d\in \cV_{c}} \Bigm( C_d(z) ( e^{ -  \langle I,\alpha_d \rangle} -1
)+1\Bigm),
\end{equation}
\begin{equation}
\label{eq:defFC} F_C ((c,p_i),\alpha,z) = p_i C_c (z) \left(1 -
\prod_{ d\in \cV_{c}} \Bigm(  C_d(z)( e^{ -  \langle I
,\alpha_d\rangle} -1 )+1\Bigm)\right).
\end{equation}
The convergence of $F_S^N$ (resp. $F_C^N$) to $F_s$ (resp. $F_C$)is
uniform in $\alpha$ and $z$, so that Assumption A1 is satisfied.It
is also easy to check that the functions $F_S$ and $F_C$ are
uniformly Lipschitz, which ensures Assumption A2.

\paragraph{Transitions of the background process $Z^N$}  Assume that the system is in state
$((c_i,p_i)_{i=1,\ldots N},\alpha,z)$. The transition kernel
$K_{\alpha}^N$ for $Z^N$ is given by: for all $z,z'\in \cZ$,

\begin{equation}\label{eq:zkernel2}
K_\alpha^N(z,z')=K_{\alpha,A_1}^N(z,z')K_{\alpha,A_2}^N(z,z')K_{\alpha,D_1}^N(z,z')K_{\alpha,D_2}^N(z,z')K_{\alpha,0}^N(z,z').
\end{equation}
$K_{\alpha,A_1}^N$, respectively $K_{\alpha,A_2}^N$, corresponds to
the transitions of links starting successful transmissions,
respectively collisions:
\begin{eqnarray*}
K_{\alpha,A_1}^N(z,z') & = & \prod_{\{c\}\in
A_1(z,z')}C_c(z)\sum_{i\in
c}{p_i\over N} \prod_{j\neq i,j\in c}(1-{p_j\over N})),\\
& = & \prod_{\{c\}\in A_1(z,z')} C_c(z)\sum_{i\in
c}\frac{p_i/N}{1-p_i/N}e^{
\langle N\log(1-\frac{I}{N}),\alpha_c\rangle})\\
&=&\prod_{\{c\}\in A_1(z,z')}C_c(z) \langle
\frac{I}{1-I/N},\alpha_c\rangle e^{ \langle
N\log(1-\frac{I}{N}),\alpha_c\rangle}),
\end{eqnarray*}

$$
K_{\alpha,A_2}^N(z,z') = \prod_{E\in A_2(z,z')}k_\alpha(e),
$$
where, if $E=\{c\}$,
\begin{eqnarray*}
k_\alpha(E) & = & C_c(z)(1-\prod_{i\in c}(1-{p_i\over N})-\sum_{i\in
c}{p_i\over N} \prod_{j\neq i,j\in c}(1-{p_j\over N})))\\
& = &C_c(z)\left(1-(1+\langle \frac{I}{1-I/N},\alpha_c\rangle) e^{
\langle N\log(1-\frac{I}{N}),\alpha_c\rangle}\right),
\end{eqnarray*}
and if $\vert E\vert \ge 2$,
\begin{eqnarray*}
k_\alpha(E) & = &\prod_{E\in A_2(z,z')} \prod_{c\in
E}C_c(z)\left(1-\prod_{i\in c}(1-{p_i\over N})\right)\\
& = &\prod_{E\in A_2(z,z')} \prod_{c\in E}C_c(z)\left(1-e^{ \langle
N\log(1-\frac{I}{N}),\alpha_c\rangle}\right).
\end{eqnarray*}
$K_{\alpha,D_1}^N$, respectively $K_{\alpha,D_2}^N$, corresponds to
the transitions of links with successful transmissions, respectively
with collisions, which become inactive:
$$
K_{\alpha,D_1}^N(z,z')=\ind_{D_1(z,z')\subset N_1(z)}\left({1\over
L}\right)^{\vert D_1(z,z')\vert},
$$
$$
K_{\alpha,D_2}^N(z,z')=\ind_{D_2(z,z')\subset N_2(z)}\left({1\over
L_c}\right)^{\vert D_2(z,z')\vert}.
$$
Finally, $K_{\alpha,0}$ corresponds to classes that are not changing
their state between $z$ and $z'$:
\begin{eqnarray*}
K_{\alpha,0}^N(z,z')&=&\left( 1-{1\over L}\right)^{\vert
N_1(z)\setminus D_1(z,z')\vert}\left( 1-{1\over L_c}\right)^{\vert
N_2(z)\setminus
D_2(z,z')\vert}\\
&&\times \prod_{c:z_c=0=z_c'}\left( C_c(z)\prod_{i\in c}(1-{p_i\over
N})+1-C_c(z)\right)\\
&=&\left( 1-{1\over L}\right)^{\vert N_1(z)\setminus
D_1(z,z')\vert}\left( 1-{1\over L_c}\right)^{\vert N_2(z)\setminus
D_2(z,z')\vert}\\
&&\times \prod_{c: z_c=0=z_c'}(C_c(z) e^{ \langle
N\log(1-\frac{I}{N}),\alpha_c\rangle}+1-C_c(z)).
\end{eqnarray*}

The limit kernel of $Z^N$ is obtained replacing $\langle
N\log(1-\frac{I}{N}),\alpha_c\rangle$ by $-\langle
I,\alpha_c\rangle$ in the above expressions. The Assumptions A3-A6 can then be easily verified.

\paragraph{Mean field asymptotics} We now verify that Assumptions A11-A12 are satisfied, implying
that Assumption A7 also holds. Let us
build a transition kernel $K$, corresponding to a process $Z$ with
values in $\cZ=\{0,1,2\}^{\vert \cC\vert}$. When equal to 0, a
component of $Z$ almost surely becomes 1 at the next slot, and
whatever the state of the system is.  The kernel $K$ then corresponds to a system
where there are always users of each class attempting to use the
channel at each slot. One can easily verify that Assumption A11-A12
are satisfied for this kernel $K$, for the partial order $\preceq$
on $\cZ$ defined by $z\preceq z'$ if and only if there is no class
$c$ such that $z_c=1$ or 2 and $z_c'=0$.

We rescale time and define $q^N_i(t)=p_i^N([Nt])$. Since the set of transitions is finite, the tightness of
$\cL (q^N_1 (\cdot))$ follows easily from Theorem 7.2 in
Ethier-Kurtz \cite{ethier} p 128. (see the comment after A8). It follows that Theorem \ref{main} applies. Assume that the class of the particle $i$ is a
r.v. fixed at the time 0 such that the vector $(c_1,\cdots,c_N)$ is
an exchangeable random vector (for example the $c_i$'s may be i.i.d.
and equal to $c$ with probability $\mu_c$). Theorem \ref{main}
asserts that as
$N\to\infty$, the $q^N_i$'s become independent and evolve according
to a measure $Q = (Q(t))_{t \in \R^+}$.

\subsection{Stationary throughputs}

\label{subsub:stat}

Assume that Assumptions A9-A10 hold, so that Theorem \ref{theost}
applies. These assumptions will be partly justified below for the
case of the binary exponential back-off algorithm. We are interested in
deriving the stationary throughputs achieved by users of various
classes. To do so, we derive the stationary distribution $Q_{\st}$
and $\pi_{Q_{\st}}$ of the particles and the background process. To
simplify the notation we write $Q_{\st}=Q$ and $\pi_Q=\pi$. Also
denote $Q_c^p=Q(\{c,p\})$ the stationary proportion of users of
class $c$ transmitting with probability $p$.

Consider the
point process of returns to the set $\cA=\{z:C_c(z)=1\}$.  Let $T_1$
denote the first return time after time zero. By the cycle formula
(see (1.3.2) in \cite{baccelli}) we may express the steady state
probability of a user in $c$ successfully transmitting a packet by
the mean time spent in the transmission state per cycle divided by
the mean cycle length.  The expectation is calculated with respect
to the Palm measure of the point process of returns to $\cA$ but in
this Markovian case this just means starting on $\cA$ with
probability $\pi^{\cA}$ which is $\pi$ renormalized to be a
probability on $\cA$.

A user in $c$ can only go into a successful transmission state once
per cycle; i.e. no other user in $c$ transmits and other users in
$\cV_c$ are either blocked or remain silent. Hence the mean time per
cycle spent in a transmission state is $\sum_{z\in
\cA}\pi^{\cA}(z)Lg(z)$ where
$$g(z)=\rho_c
\prod_{d\in\cV_c, d\neq c}(C_d(z)(e^{-\rho_d}-1)+1).$$
Moreover, $\sum_{z\in \cA}\pi^{\cA}(z)E_z[T_1]=\frac{1}{\pi({\cal A})}$;
i.e. the intensity of the point process of visits to $\cA$. Finally
 the total throughput of the users of
class $c$ is
\begin{equation}\label{eq:thru}
\gamma_c=\sum_{z:C_c(z)=1}\pi(z)L \rho_c\prod_{d\in
\cV_c}\left(C_d(z)(e^{-\rho_d}-1)+1\right),
\end{equation} where
\begin{equation}\label{eq:rho}
\rho_c=\sum_{p\in\cB}pQ_c^p,
\end{equation}
which can be interpreted as the probability that a user of class $c$
attempts to use the channel at the end of an empty slot. We now
evaluate $Q$ and $\pi$. Note that $\pi$ depends on $Q$ through the
$\rho_c$'s only (see (\ref{eq:zkernel2}) and its limiting
expression). Then we can write:
\begin{equation}\label{eq:pi}
\pi(z)=\Phi(z,\rho_c,c\in\cC); \mbox{i.e. $\pi$ is a function of $\rho_c,c\in\cC$}, 
\end{equation}

Now define $G_c$, $H_c$ and $I_c$ as follows:
\begin{eqnarray}
G_c &=& \sum_z \pi(z)  C_c (z) \prod_{ d\in \cV_{c}} \Bigm( C_d(z) (
e^{ - \rho_d} -1 )+1\Bigm),\label{eq:fin1}\\
H_c &=& \sum_z \pi(z) C_c (z) \left(1 - \prod_{ d\in \cV_{c}} \Bigm(
C_d(z)( e^{ - \rho_d } -1 )+1\Bigm)\right),\label{eq:fin2}\\
I_c &=& G_c+H_c=\sum_z \pi(z)C_c(z).
\end{eqnarray}
$G_c, H_c, I_c$ depend on $Q$ through the $\rho_c$'s only. We have
for all $c,p$: $pG_c = \overline F_S ( (c,p) ,Q )$, $pH_c= \overline
F_C( (p,c),Q)$. The marginals $Q^p_c$ satisfy the balance equations
(\ref{eq:balance}), i.e., for all $c,p$,
\begin{equation}\label{eq:qq}
G_c\left(\sum_{p'\in\cB: S(p')=p}p'Q_c^{p'}-
pQ_c^p\right)+H_c\left(\sum_{p'\in\cB: C(p')=p}p'Q_c^{p'}-
pQ_c^p\right)=0.
\end{equation}
They also satisfy:
\begin{equation}\label{eq:qq2}
\forall c\in\cC,\quad \sum_{p\in\cB} Q_c^p = \mu_c.
\end{equation}
Summarizing the above analysis, we have:
\begin{theorem}\label{th:final}
The stationary distribution $Q$ is characterized by the set of
equations (\ref{eq:rho}), (\ref{eq:pi}), (\ref{eq:fin1}),
(\ref{eq:fin2}), (\ref{eq:qq}), (\ref{eq:qq2}).
\end{theorem}

\subsection{The binary exponential back-off algorithm}
We now examine the specific case of the binary exponential back-off
algorithm. We first justify Assumption A9.

\subsubsection{Tightness of stationary distributions}

\begin{lemma}\label{lem:rec} In case of the exponential back-off algorithm, there exists a $p^*>0$, such that for any $0< p_0 < p^*$, the
Markov process $(X_i^N(k),Z^N(k))_{k\in \N}$ is positive recurrent
for all $N$ and the family of stationary distributions $\cL_{st}(X_1^N(0))$
is tight.
\end{lemma}

Deriving a tight bound for $p^*$ would involve technical details
which are beyond the scope of this paper.We will only sketch the
main idea and prove $p^* >0$. Also to clarify the presentation, we
assume here that $L=L_c$. Along the proof of Lemma \ref{lem:rec}, we
may check that the statement of Lemma \ref{lem:rec} holds for $p^* =
\frac{\ln 2 }{ L \overline \mu }$, where $\overline \mu= \max_{c \in
\cC}  \overline \mu_c $ and $\overline \mu_c = \sum_{d\in \cV_c}
\mu_d$ is the mean proportion of particles which are in interaction
with particles of class $c$.

\bp To prove the recurrence we introduce a fictive system which  stochastically bounds  $p_1^N(k)$.

In the fictive system, the states of the particles $i \geq 2$ are independent, a particle $i \geq 2$ has two states: active or inactive. If the particle $ i \geq 2$, is active, it remains active for the next slot with probability $1-1/L$, if it is inactive, it becomes active with probability $p_0 /N$.  The stationary probability that the particle $i$ is active is $L / (L + N /p_0)$ and the stationary probability that at least one is active is $a_N = 1 - (1 - L/ ( L + N/p_0) )^{N-1} $ which converges to $a = 1 - e^{-Lp_0}$.

The particle $1$ tries to become active at slot $k$ with probability $ p_1^N(k)/N$. If it remains inactive, $p^N_1 (k)  = p^N_1 (k+1)$. If it is active and if another  particle is also active, then the particle $1$ encounters a collision and $p_1^N (k+1) = p_1 ^ N (k) /2$. Otherwise $p_1 ^ N (k+1) = p_0$.

Clearly, this virtual system is stochastically less than or equal to
$p_1^N(k)$ in the exponential back-off case.

Let $b^N (k) = p_0 / p_1 ^N (k)$, $b^N (k) \in \{ 2^n \}_{ n \in \N}$, the lemma will follow if we prove that for $p_0$ small enough,
\begin{equation}
\label{eq:bN}
\sup_{N,k} \EE [  b ^N (k)   \, | \, b ^ N (0) = 1 ] < \infty.
\end{equation}
In the remaining part of the proof, using elements of queueing theory, we justify (\ref{eq:bN}).

We first analyze the sequence of slots such that none of the particles $i \geq 2$ is active.
If the particle $i \geq 2$ is active at time $k$, let $l_i (k)$ be the number of slots the particle remains active. $l_i(k)$ is a geometric distribution with parameter $1/L$.  Now, let
$$W^N(k) = \max_{2 \leq i \leq N}  \chi\{ i \hbox{ active} \} l_i (k).$$
If $W^N (k) = 0$ none of the particles $i \geq 2$ is active at time $k$.  $W^N$ satisfies the recursion:
$$
W^N(k+1)  =  \max \Bigm( W^N(k) -1 , \max_{2 \leq i \leq N} \chi\{\hbox{$i$ active at $k+1$, inactive at $k$}\} l_i (k+1) \Bigm) .
$$
$ W^N$ is thus the workload in a $G/G/\infty$ queue with inter-arrival time $1$ and service time requirement $\sigma^N (k+1) = \max_{i \geq 2} \chi\{\hbox{$i$ active at $k+1$, inactive at $k$}\} l_i (k+1)$. Independently of the past, $\sigma^N(k+1)$  is easily bounded stochastically; indeed,  let $0 < s < \ln L$,
\begin{eqnarray*}
\EE e^{s \sigma^N(k+1)} & \leq  & 1 + \sum_{i=2}^N \EE  \chi\{\hbox{$i$ active at $k+1$, inactive at $k$}\} e^{s l_i(k+1)} \\
& \leq & 1 + (N-1) \frac{p_0}{N} \EE e^{s l_i (k+1)} \\
& \leq & 1+p_0\frac{e^s/L}{1-(1/L)e^s}
\end{eqnarray*}
Note that this last bound is uniform in $N$ and $k$. Let $\theta_0 = 0 $, $\theta_{n+1} = \inf \{ k > \theta_{n} : W^N(k) = 0\}$, and $\Theta^N = \{ \theta_n\}_{n \in \N}$.  Classically, there exists $C>0$ such that for all $N$:
$$
\EE[ e^{ C (\theta_{n+1} - \theta_n) }  \, | \, W^N(0) = 0 ] < \infty,
$$
see for example Appendix A.4 in \cite{rstAAP}.  By the renewal theorem, we deduce,  uniformly in $N$,  $\lim_{k \to \infty} \PP ( k \in \Theta^N) = \frac{1}{\EE \theta_1} = 1 - a_N$. Moreover, the monotonicity of $W^N (k)$ with respect to the initial condition implies easily that $\PP ( k \in \Theta^N | W^N (0) = 0 ) \geq \lim_{k \to \infty} \PP ( k \in \Theta^N)   =1 - a_N$. Since $1 - a_N$ converges to $e^{-Lp_0}$, it follows that
\begin{equation}
\label{eq:lambdap0}
\lim_{p_0 \to 0} \inf_{k,N} \PP ( k \in \Theta^N | W^N(0) = 0 ) = 1.
\end{equation}
We now turn back to the process $b^N$ and prove (\ref{eq:bN}). Let $U(k)$ be a sequence of independent and uniformly distributed variables on $[0,1]$. We may write
\begin{eqnarray*}
b^N(k+1) & = & b^N(k) \chi_{\{ U(k+1) > \frac{p_0}{ b^N (k) N} \} } + 2 b^N(k)   \chi_{\{ U(k+1) \leq \frac{p_0}{ b^N (k) N} \}} \chi_{\{ k \notin \Theta^N\} }\\
& & \quad \quad \quad \quad \quad  +   \chi_{\{U(k+1) \leq \frac{p_0}{ b^N (k) N} \}} \chi_{\{k \in \Theta^N\}}.
\end{eqnarray*}
In particular
\begin{eqnarray*}
b^N(k+1) \chi_{\{ b^N(k) \geq 2\}} & \leq & b^N(k)  \chi_{\{ U(k+1) > \frac{p_0}{ b^N (k) N} \} }+ 2 b^N(k)    \chi_{\{ U(k+1) \leq \frac{p_0}{ b^N (k) N} \} }  \chi_{\{ k \notin \Theta^N\}} \\
& & \quad \quad \quad \quad \quad +    \chi_{\{ U(k+1) \leq \frac{p_0}{ 2 N} \}}  \chi_{\{ k \in \Theta^N\}}
\end{eqnarray*}
Taking expectation, we obtain
\begin{eqnarray*}
\EE b^N(k+1) \chi_{\{b^N(k) \geq 2\}}   & \leq &  \EE b^N(k)- \frac{p_0}{N}  + 2 \frac{p_0}{N}  \PP ( k \notin \Theta^N)  + \frac{p_0}{2N} \PP( k \in \Theta^N)    \\
& \leq &  \EE b^N(k) - \frac{p_0}{N} \left( \frac{3}{2} \PP ( k \in \Theta^N) -1\right) .
\end{eqnarray*}
Similarly, since $b^N(k+1) \chi_{\{b^N (k) = 1\} } \leq 2$, we have:
$$
\EE b^N(k+1) \leq   \max\left(2,\EE b^N(k) - \frac{p_0}{N} \left( \frac{3}{2} \PP ( k \in \Theta^N) -1\right)\right).
$$
From (\ref{eq:lambdap0}), for $p_0$ small enough, for all $N$ and $k
\geq 0$, $\PP ( k \in \Theta^N) > 2/3$. We deduce by recursion that
$ \EE [ b^N (k) | b^N(0) = 1 ]  \leq 2$ and (\ref{eq:bN}) holds. \ep

\subsubsection{Stationary distribution}

Now Lemma \ref{lem:rec} implies that Assumption A9 holds. It remains to check Assumption A10. In the next paragraph, we state that Assumption A10 holds if 
there is a unique class of users, i.e., in the case of full interference. For the general case of partial interference, we can only provide a
characterization of the equilibrium point of the dynamical system
(\ref{eqdiff1}). We leave the study of its global stability for
future work.

So let us assume that an equilibrium point exists, and denote by $Q$
this point. Further define $Q_c^n=Q(\{c,p_02^{-n}\})$ for all
$n\in\N$. Then we have:
\begin{equation*}
Q^{n-1}_c \overline F_C ( ( p_02^{-n+1} ,c),Q )= Q^{n}_c ( \overline
F_S ( (p_0 2^{-n},c),Q)+ \overline F_C ( (p_0 2^{-n},c),Q)),
\end{equation*}
or equivalently
\begin{equation}\label{simple}
2Q^{n-1}_c H_c=Q^{n}_c I_c,
\end{equation}
and
\begin{equation*} \sum_{n\geq 0}Q^{n}_c \overline F_S (
(p_02^{-n},c) ,Q)=Q^{0}_c \overline  F_C ( (p_0,c),Q)),
\end{equation*}
or equivalently
\begin{equation}\label{simple0}
\rho_c G_c =Q^{0}_c H_c.
\end{equation}

Solving (\ref{simple}) and (\ref{simple0}) leads to a solution of
the form $Q^{n}_c=\beta_c (2H_c/I_c)^{n}$. Since $\sum_n
Q^n_c=\mu_c$, we have $\beta_c=Q^0_c = \mu_c ( 1 - 2 H_c / I_c ) \mu_c (G_c-H_c) / I_c$. We require that $H_c<G_c$ or equivalently,
$G_c/I_c>1/2$. $G_c/I_c$ may be interpreted as the probability in
steady state that no user of class in $\cV_c$ tries to access the
channel given that no user of class in $\cV_c$ are currently
sending. Next, $\rho_c=\sum_{n\geq 0}p_02^{-n}Q^n_c$, which implies
that:
\begin{equation}\label{big}
\rho_c=p_0\mu_c\frac{G_c -  H_c}{G_c} .
\end{equation}
Now the following corollary summarizes the above analysis and then
it characterizes the system behavior in steady state and in case of
exponential back-off algorithms.
\begin{coro}
\label{cor:statQ}
The stationary distribution $Q$ is given by: for all $c\in \cC$,
$$
Q_c ^ n = \mu_c {G_c-H_c\over G_c+Hc}\left({2H_c\over
G_c+H_c}\right)^n,
$$
where the $G_c's$, $H_s$'s, and $\rho_c$'s are the unique solutions
of the system of equations  (\ref{eq:pi}), (\ref{eq:fin1}),
(\ref{eq:fin2}), (\ref{big}).
\end{coro}

\subsubsection{Global stability in the mean field regime for networks with full interference}\label{globalstab}

In this paragraph, we consider the exponential backoff algorithm and we assume moreover that there is a unique class of users. In that case, the analysis in greatly simplified: the environment variable $Z^N$ and the clear-to-send function $C$ are then identical for all users. The backoffs of the users evolves only if  $Z^N = 0$. Thus, up to sampling by the times such that $Z^N(t) = 0$, in order to analyze the backoff process, we may assume without loss of generality that $L = L_c =1$.

Let $Q^n(t)$ be the mean field limit of the proportion of users with backoff $p_0/2^n$.  From Theorem \ref{main}, given an initial distribution $\{Q^n(0);
n=0,1,\ldots\}$, the limit evolves as, for all $n\geq 1$,
\begin{equation}\label{eqdiff1b}
{dQ^n\over
dt}(t)=2^{1-n}p_0 Q^{n-1}(t)\big(1-\exp(-\sum_{i=0}^{\infty}
2^{-i}p_0Q^i(t))\big)- 2^{-n} p_0 Q^n(t),
\end{equation}
\begin{equation}\label{eqdiff2b}
{dQ^0\over
dt}(t)=\sum_{n=0}^{\infty}2^{-n}p_0Q^n(t)\exp(-\sum_{i=0}^{\infty}2^{-i}p_0Q^i(t))-p_0Q^0(t).
\end{equation}

Following \S \ref{subsub:stat} and Corollary \ref{cor:statQ}, the dynamic system described by
differential equations (\ref{eqdiff1b})-(\ref{eqdiff2b}) admits a
unique equilibrium point $Q_{st}=\{Q^n_{st};n=0,1,\ldots\}$ defined by:
$$
\forall n\ge 0,\quad Q^n_{st}=(2(1-e^{-\rho_{st}}))^nQ^0_{st},\quad Q_{st}^0=\rho_{st}
e^{-\rho_{st}}/p_0
$$
where $\rho_{st}$ solves $p_0e^\rho+\rho-2p_0=0$ or
$p_0=\rho/(2-e^\rho)$.
Note that $\rho_{st} < \ln(2)$ so necessarily  $2(1-e^{-\rho_{st}})< 1$, and the stationary distribution always exists. Moreover $\rho_{st}=\sum_{i=0}^{\infty}2^{-i}p_0Q_{st}^i$.

Now let $\rho(t)=\sum_{i=0}^{\infty}2^{-i}p_0Q^i(t)$ so (\ref{eqdiff1b})-(\ref{eqdiff2b})  can
be written as
\begin{equation}\label{eqdiff1d}
\frac{dQ^n(t)}{dt}
=2^{1-n}p_0Q^{n-1}(t)\big(1-e^{-\rho(t)}\big)-2^{-n}p_0 Q^n(t),\quad \hbox{ for all
}n\ge 1,
\end{equation}
\begin{equation}\label{eqdiff2d}
\frac{dQ^0(t)}{dt}=\rho(t)\exp(-\rho(t))-p_0 Q^0(t).
\end{equation}

In complete interaction, Assumption A9 holds. Indeed, we have the following:
\begin{theorem}\label{limit}
If $p_0 < \ln(2)$, for any initial condition $Q(0)$, $Q(t)$ converges (weakly) to the measure $Q_{st}$.
\end{theorem}

\begin{lemma}\label{tight}
If $p_0 < \ln(2)$, the sequence of measures $Q(t)=\{Q^n(t);n=0,1,\ldots\}$ is tight.
\end{lemma}

\noindent {\em Proof of Lemma \ref{tight}.}
We define the linear system,
\begin{equation}\label{eqdiff1c}
{dB^n\over
dt}(t)=2^{1-n}p_0 B^{n-1}(t)\big(1-\exp(-p_0)\big)- 2^{-n} p_0 B^n(t),\; \hbox{ for all
}n\ge 1,
\end{equation}
\begin{equation}\label{eqdiff2c}
{dB^0\over
dt}(t)=\sum_{n=0}^{\infty}2^{-n}p_0B^n(t)\exp(-p_0)-p_0B^0(t)
\end{equation}
with initial condition $B^n(0)=Q^n(0)$ for all $n$. First note that the time derivative of $\sum_n B^n(t)$
is zero, hence $\sum_{n=0}^{\infty}B^n(t)=1$ for all $t \geq 0$. Note also that $\rho(t) \leq  p_0$.
$B^n(t)$ corresponds to mean field limit of the proportion of users with backoff $p_0 2^{-n}$ when each user is in interaction with $N$ other users with backoff $p_0$.
We may then check that the probability measure $B(t)=\{B^n(t);n\geq 0\}$ is stochastically larger than $Q(t)$: for all $m \geq 1$, $\sum_{n \geq m} B^n(t) \geq \sum_{n \geq m } Q^n(t)$. However  $B(t)$ converges to the unique invariant probability measure of the linear system: $B_{st}^n = (2 ( 1 - \exp( - p_0))^n B_{st}^0$ (recall that $p_0 < \ln (2)$). Since $B(t)$ converges, it is therefore tight.  It follows that $Q(t)$ is tight. \ep

\noindent {\em Proof of Theorem \ref{limit}.} Let $\liminf_{t\to\infty}\rho(t)=\rho_b$.  Pick a subsequence $t_k$
such that $\lim_{t_k\to\infty}\rho(t_k)=\rho_b$ and such that the limit $\lim_{t_k\to\infty}Q^n(t_k)=Q^n(\infty)$
exists for all $n$. By Lemma \ref{tight}, $Q(\infty)=\{Q^n(\infty);n=0,1,\ldots\}$
is a probability measure and $\sum_{i=0}^{\infty}2^{-i}p_0Q^i(\infty))=\rho_b$.

Let  $f_b(t)= \inf_{u\geq t}\rho(u)$.
Note that $f_b(t)$ increases to $\rho_b$ and $f_b(t) \leq \rho(t)$ for all $t \geq 0$.
Now consider the system
\begin{equation*}\label{eqdiff1e}
\frac{d\tilde{Q}^n(t)}{dt}
=2^{1-n}p_0\tilde{Q}^{n-1}(t)\big(1-e^{-f_b(t)}\big)-2^{-n}p_0 \tilde{Q}^n(t),\quad \hbox{ for all
}n\ge 1,
\end{equation*}
\begin{equation*}\label{eqdiff2e}
\frac{d\tilde{Q}^0(t)}{dt}=f_b(t)\exp(-f_b(t))-p_0 \tilde{Q}^0(t)
\end{equation*}
with initial condition $\tilde{Q}^n(0) = Q^n (0)$ for all $n \in \N$.
Now notice
that the function $\rho \exp(-\rho)$ is strictly increasing for $0\leq
\rho\leq 1$. Hence, for $t\geq 0$,
\begin{eqnarray}
\tilde{Q}^0(t)&=&e^{-p_0 t}\tilde{Q}^0(0)+e^{-p_0 t}\int_0^t e^{p_0 s}f_b(s )\exp(-f_b(s))ds\label{for0}\\
&\leq &e^{-p_0 t}Q^0(0)+e^{-p_0t}\int_0^t e^{p_0 s} \rho(s) \exp(-\rho(s)) ds \nonumber\\
&=&Q^0(t)  \nonumber.
\end{eqnarray}
Therefore for all $t \geq 0$,
$$
0<\tilde Q^0 (t) \leq Q^0 (t).
$$
We then prove by recursion on $n$ that
\begin{equation}
\label{eq:rectilde}
\forall n \in \N, \; \forall t \geq 0, \quad 0<\tilde Q^n (t) \leq Q^{n-1}(t).
\end{equation}
Let $n\geq 1$, and assume that for all $t \geq 0$, $\tilde Q^{n-1} (t)  \leq Q^n(t)$. We have:
\begin{eqnarray*}
\tilde{Q}^n(t)&=&e^{-p_02^{-n}t}\tilde{Q}^n(0)+e^{-p_02^{-n}t}\int_0^t e^{p_02^{-n}s}p_0 2^{1-n}\tilde{Q}^{n-1}(s)(1-e^{-f_b(s)})ds\nonumber\\
&\leq &e^{-p_02^{-n}t}Q^n(0)+e^{-p_02^{-n}t}\int_0^t e^{p_02^{-n}s}p_0 2^{1-n}Q^{n-1}(s)(1-e^{-\rho(s)})ds\nonumber\\
&=&Q^n(t).\label{forn}
\end{eqnarray*}
From (\ref{forn}) we also conclude that $\tilde{Q}^n(t)\leq Q^n(t)$ for all $t$ so (\ref{eq:rectilde}) follows. Next, using L'H\^{o}pital's rule in (\ref{for0}) we get
\begin{eqnarray*}
\tilde{Q}^0(\infty)&:=&\lim_{t\to\infty}\tilde{Q}^0(t)\\
&=&\lim_{t\to\infty}\frac{ e^{p_0 t}f_b(t)\exp(-f_b(t)) }{ p_0e^{p_0 t} }=\frac{\rho_b e^{-\rho_b} }{p_0}.
\end{eqnarray*}
Moreover, by iteration
\begin{eqnarray*}
\tilde{Q}^n(\infty)&:=&\lim_{t\to\infty}\tilde{Q}^n(t)\\
&=&\lim_{t\to\infty}\frac{e^{p_0 2^{-n}t}p_0 2^{1-n}\tilde{Q}^{n-1}(t)(1-e^{-f_b(t)})}{p_02^{-n}e^{p_0 2^{-n}t}}\\
&=&2(1-e^{-\rho_b})\lim_{t\to\infty}\tilde{Q}^{n-1}(t)=(2(1-e^{-\rho_b}))^n\frac{\rho_b e^{-\rho_b}}{p_0}.
\end{eqnarray*}
This expression of $\tilde{Q}(\infty)$ implies
$$\sum_{n=0}^{\infty}p_02^{-n}\tilde{Q}^n(\infty)=\rho_b=\sum_{n=0}^{\infty}p_02^{-n}Q^n(\infty).$$
However, (\ref{eq:rectilde}) implies that  for all $n \in \N$, $\tilde{Q}^n(\infty)\leq Q^n(\infty)$, it follows that  $\tilde Q^n(\infty) = Q^n (\infty)$, and therefore,
$$
\sum_{n=0}^{\infty} \tilde Q^n(\infty)=1.
$$
But from the above expression for $\tilde{Q}^n(\infty)$, we have:
\begin{eqnarray*}
\sum_{n=0}^{\infty}\tilde{Q}^n(\infty)&=&\frac{\rho_b e^{-\rho_b}}{p_0}\frac{1}{1-2(1-e^{-\rho_b})}\\
&=&\frac{\rho_b}{p_0(2-e^{\rho_b})}.
\end{eqnarray*}
We conclude $\rho_b$ solves $p_0  = \rho / (2 - e^{\rho})$. Hence $\rho_b= \rho_{st}$ and this means
$$Q^n(\infty)=\tilde{Q}^n(\infty)=Q^n_{st}\quad \mbox{ for all $n$. }$$
Hence any subsequence of $Q^n(t)$ converges to $Q^n_{st}$ and this gives our result. \ep

\subsection{A numerical example}\label{sec:examplenum}

\begin{figure}[h]
\centering
\includegraphics[width=8.5cm]{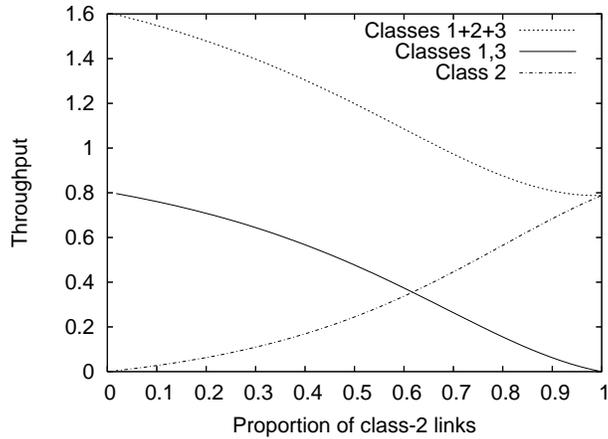}
\caption{\label{fig:results1}Throughputs of links as a function of
$\mu_2$ - $L=L_c=100$ slots.}
\end{figure}

\begin{figure}[h]
\centering
\includegraphics[width=8.5cm]{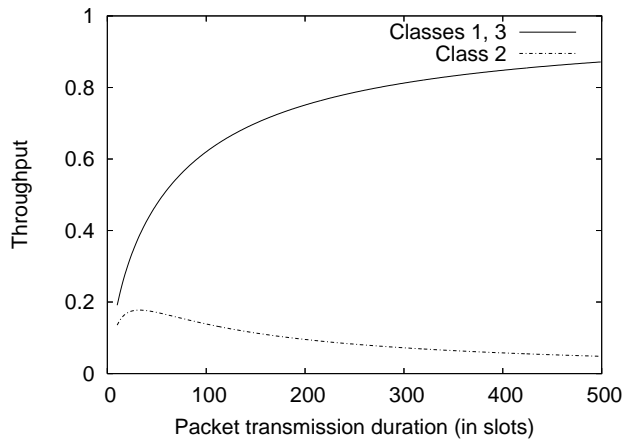}
\caption{\label{fig:results2}Throughputs of links as a function of
$L$, when $\mu_1=\mu_2=\mu_3$.}
\end{figure}

We now illustrate our analytical results on the simple network of Figure \ref{firstpic}. Each link runs an exponential back-off algorithm with
$p_0=1/16$ as specified in the 802.11 standard \cite{ieee}. In Figure
\ref{fig:results1}, the throughputs of the various user classes are
presented assuming that the proportions of users of class 1 and 3
are identical, $\mu_1=\mu_3$. We assume here that $L_c=L$. We give
the throughputs as a function a the proportion of users of class 2.
Here the packet duration is fixed and equal to $L=100$ slots. The
total network throughput decreases when the proportion of class-2
users increases, which illustrates the loss of efficiency due to the
network spatial heterogeneity. In Figure \ref{fig:results2}, we
assume a uniform user distribution among the 3 classes,
$\mu_1=\mu_2=\mu_3$, and we give the throughputs as a function of
the packet duration $L$. First note that whatever the value of $L$,
the network is highly unfair: for example when $L=100$ slots, the
throughput of a user of class 1 is almost 5 times greater than that
of a user of class 2. This unfairness increases with $L$ and
ultimately when $L$ is very large, users of class 2 never access the
channel successfully. We have verified through simulation that mean field asymptotics led to quite accurate performance approximations, even in the case of systems with a small number of users. This has been also observed in \cite{bianchi} for networks with full interference.

\section*{Acknowledgments}
This work was done during C. Bordenave's and A. Proutiere's visit to University of Ottawa. We wish to thank Prof. McDonald for his kind hospitality.


\begin{thebibliography}{9}

\bibitem{ethernet}
IEEE Standard for Information technology � Telecommunications and
information exchange between systems � Local and metropolitan area
networks, IEEE std., 1983.


\bibitem{ieee}
IEEE Standard for Wireless LAN Medium Access Protocol and Physical
Layer Specifications, IEEE std., August 1999.

\bibitem{abramson} \textsc{Abramson, N.} (1970). The ALOHA system - Another alternative for computer communications, in: {\it Proc. of AFIPS press} \textbf{37}.

%\bibitem{aldous} \textsc{Aldous, D.} (1986). Ultimate instability of exponential back-off protocol for acknowledgement-based transmission control of random
%access communication channels, {\it IEEE trans. on Information
%Theory} \textbf{33-2}.


\bibitem{rstAAP}
\textsc{Baccelli, F.} and \textsc{Bordenave, C.} (2007). The radial spanning tree of a Poisson point process {\it 
Ann. Appl. Probab.} Vol. 17, No. 1, 305-359.



\bibitem{baccelli}
\textsc{Baccelli, F.} and \textsc{Br\'{e}maud, P.} (1994). {\it
Elements of Queueing Theory: Palm-Martingale Calculus and Stochastic
Recurrences}, Springer-Verlag, New York.

\bibitem{ben} \textsc{Bensoussan, A.,}  \textsc{Lions, J-L.} and
\textsc{Papanicolaou, G.}(1978).
Asymptotic analysis for periodic structures.
{\it Studies in Mathematics and its Applications}, {\bf 5}. North-Holland
Publishing Co., Amsterdam-New York.

\bibitem{bianchi} \textsc{Bianchi, G.} (2000). Performance Analysis of the IEEE 802.11
Distributed Coordination Function, {\it IEEE Journal on Selected
Areas in Communications} \textbf{18-3} 535-547.

\bibitem{billingsley}
\textsc{Billingsley, P.} (1999). \textit{Convergence of Probability
Measures}, 2nd ed. Wiley, New York.

\bibitem{allerton}
\textsc{Bordenave, C.} and \textsc{McDonald, D.} and
\textsc{Proutiere, A.} (2005). Random multi-access algorithms: A
mean field analysis, In: {\it Proc. of the 43th Allerton conference
on communication, control, and computing.}

\bibitem{ussquared}\textsc{Bordenave, C.} and \textsc{McDonald, D.} and
\textsc{Proutiere, A.} (2008). Asymptotic stability region of slotted-Aloha. Submitted
and archived:arXiv:0809.5023.

\bibitem{cordier}
\textsc{Cordier, S.} and \textsc{Pareschi, L.} and \textsc{Toscani, G.} (2005). On
a kinetic model for a simple market economy, {\it Journal of
Statistical Physics} {\bf 120}.


\bibitem{dawson}
\textsc{Dawson, D. A.} (1983). Critical Dynamics and Fluctuations for a Mean-Field model of cooperative
behavior, {\it Journal of Statistical Physics} {\bf 31} 29-85.

\bibitem{donnelly}
\textsc{Donnelly, P.} and \textsc{Kurtz, T.} (1999).
Particle representations for measure-valued population models, {\it
Annals of Probability} {\bf 27} 166-205.

\bibitem{DT06}
\textsc{Durvy, M.} and \textsc{Thiran, P.} (2006). A Packing
Approach to Compare Slotted and Non-slotted Medium Access Control,
In: {\it proc. of IEEE Infocom 2006}.

\bibitem{ethier}
\textsc{Ethier, S.} and \textsc{Kurtz, T.} (1985). {\it Markov
processes}, Wiley, New york.

\bibitem{feller}
\textsc{Feller, W.} (1971). {\it Introduction to Probability Theory and its Applications}, vol 2, Wiley, New York.



\bibitem{graham97}
\textsc{Graham, C.} and \textsc{M\'el\'eard, S.} (1997). Stochastic
particle approximations for generalized Boltzmann models and
convergence estimates, {\it Annals of Probability} \textbf{28}
115-132.

\bibitem{graham}
\textsc{Graham, C.} (2000). Chaoticity on path space for a queueing
network with selection of the shortest queue among several, {\it
Journal of Applied Probability} \textbf{37} 198-211.

%\bibitem{hajek}
%\textsc{Hajek, B.} and \textsc{van Loon, T.} (1982). Decentralized
%Dynamic Control of a Multiaccess Broadcast Channel, {\it IEEE trans.
%on Automatic Control} \textbf{27-3}.

%\bibitem{hastad}
%\textsc{H{\aa}stad, J.} and \textsc{Leighton, T.} and
%\textsc{Rogoff, B.} (1996). Analysis of backoff protocols for
%multiple access channels, {\it SIAM Journal of Computing}
%\textbf{25-4}.

\bibitem{jiang}
\textsc{Jiang, L.} and \textsc{Liew, S.C.} (2008). Improving Throughput and Fairness by Reducing Exposed and Hidden Nodes in 802.11 Networks, {\it IEEE transaction on Mobile Computing}, January.




\bibitem{karatzas}
\textsc{Karatzas, Y.} and \textsc{Shubik, M.} and \textsc{Sudderth,
W.} (1997). A Strategic Market Game with Secured Lending, {\it
Journal of Mathematical Economics} \textbf{28} 207-247.



\bibitem{kelly}
\textsc{Kelly, F.P.} (1985). Stochastic models of computer
communication systems, {\it Journal of the Royal Statistical
Society} \textbf{47-3}.


\bibitem{kellyloss}
\textsc{Kelly, F.P.} (1991). Loss Networks, {\it The Annals of
Applied Probability} \textbf{1-3} 319-378.

\bibitem{kurtzaverage}
\textsc{Kurtz, T. G.} (1992).  Averaging for martingale problems and stochastic approximation, {\it Lecture Notes in Control and Information Sciences} \textbf{177},
Springer Berlin / Heidelberg.


\bibitem{meyntweedie}
\textsc{Meyn, S.} and \textsc{Tweedie, R.} (1993). {\it Markov Chains and Stochastic Stability}, Springer-Verlag,
New York.

\bibitem{MS02}
\textsc{M\"uller, A.} and \textsc{Stoyan, D.} (2002). {\it
Comparison Methods for Stochastic Models and Risks}, Wiley, New
York.

\bibitem{shiga}
\textsc{Shiga, T.} and \textsc{Tanaka, H.} (1985). Central limit
theorem for a system of Markovian particles with mean field
interaction, {\it Z. Wahrscheinlichkeitsth.} \textbf{69} 439-459.

\bibitem{sznitman}
\textsc{Sznitman, A.S.} (1991). Propagation of chaos, in: {\it Ecole
d'\'et\'e de probabilit\'es de Saint-Flour XIX}, lecture notes in
Maths 1464, Springer, Berlin.

\end{thebibliography}
\end{document}